\documentclass[final]{autart}

\usepackage{amsmath}
\usepackage{amssymb}
\usepackage{makeidx}
\usepackage{multicol}
\usepackage[bottom]{footmisc}
\usepackage{graphicx}
\usepackage{color}
\usepackage{natbib}
\usepackage{epstopdf}
\usepackage{algorithm}
\usepackage{algcompatible}
\usepackage{calc}

\newtheorem{approxi}{Approximation}
\newtheorem{assump}{Assumption}

\newcommand{\app}{\begin{approxi}}
\newcommand{\eapp}{\end{approxi}}
\newcommand{\ass}{\begin{assump}}
\newcommand{\eass}{\end{assump}}
\newcommand{\teo}{\begin{thm}}
\newcommand{\eteo}{\end{thm}}
\newcommand{\corr}{\begin{cor}}
\newcommand{\ecorr}{\end{cor}}
\newcommand{\pro}{\begin{prop}}
\newcommand{\epro}{\end{prop}}
\newcommand{\lemma}{\begin{lem}}
\newcommand{\elemma}{\end{lem}}
\newcommand{\pb}{\begin{prob}}
\newcommand{\epb}{\end{prob}}
\newcommand{\df}{\begin{defn}}
\newcommand{\edf}{\end{defn}}
\newcommand{\rema}{\begin{rem}}
\newcommand{\erema}{\end{rem}}

\newcommand{\al}[1]{\begin{align} #1 \end{align}}
\newcommand{\nn}{\nonumber}

\newcommand{\tr}{\mathop{\rm tr}}  

\newcommand{\Ec}{ \mathcal{E}}

\newcommand{\Hc}{ \mathcal{H}}
\newcommand{\Ic}{ \mathcal{I}}

\newcommand{\Lc}{ \mathcal{L}}
\newcommand{\Mc}{ \mathcal{M}}

\newcommand{\Pc}{ \mathcal{P}}

\newcommand{\Sc}{ \mathcal{S}}

\newcommand{\Cs}{ \mathbb{C}}

\newcommand{\Es}{ \mathbb{E}}

\newcommand{\Ns}{ \mathbb{N}}

\newcommand{\Rs}{ \mathbb{R}}

\newcommand{\Zs}{ \mathbb{Z}}

\newcommand{\yv}{\mathrm{y}}
\newcommand{\aaa}{\mathrm{a}}
\newcommand{\bbb}{\mathrm{b}}
\newcommand{\ccc}{\mathrm{c}}
\newcommand{\ev}{\mathrm{e}}

\newcommand{\pp}{\mathbf{p}}

\begin{document}

\begin{frontmatter}
\title{Sparse plus Low rank Network Identification:\\ A Nonparametric  Approach\thanksref{footnoteinfo}}

\thanks[footnoteinfo]{This work has been partially supported by the FIRB project ``Learning
meets time'' (RBFR12M3AC) funded by MIUR.}

\author[Padova]{Mattia Zorzi}\ead{zorzimat@dei.unipd.it},
\author[Padova]{Alessandro Chiuso}\ead{chiuso@dei.unipd.it}

\address[Padova]{Dipartimento di Ingegneria dell'Informazione, Universit\`a degli studi di
Padova, via Gradenigo 6/B, 35131 Padova, Italy}


\begin{abstract}
Modeling and identification of high-dimensional stochastic processes is ubiquitous in many fields. In particular, there is a growing interest
in modeling stochastic processes with simple and interpretable structures. In many applications, such as econometrics and biomedical sciences, it seems natural to describe each component of that stochastic process in terms of few factor  variables, which are not accessible for observation, and possibly of few other components of the stochastic process. These relations can be encoded in graphical way via a structured dynamic network, referred to as  ``sparse plus low-rank (S+L) network'' hereafter. The problem of finding the S+L network as well as the dynamic model can be posed as a system identification problem. 
 In this paper, we introduce two new nonparametric methods to identify dynamic models for stochastic processes described by a  S+L network.  These methods take inspiration from   regularized  estimators based on recently introduced kernels (e.g. ``stable spline'',  ``tuned-correlated''  etc.). Numerical examples show the benefit to introduce the S+L structure in the identification procedure. 
\end{abstract}
\end{frontmatter}

\section{Introduction}

In many applications, high-dimensional  data are measured to describe the underlying  phenomena. These data are typically measured over time and can thus be modeled as a stochastic process whose components are referred to as \emph{manifest} variables, i.e.  accessible for observation. Very often, in high dimensional time series modeling it is necessary to control the model complexity to be able to obtain sensible results from a finite set of measured data. In addition, models should be interpretable in the sense of providing an insight into the data generation mechanism.  

One possible way to limit the complexity is to postulate, often very reasonably,  that these observations share some common behaviour (comovements) which, in turn, can be described by a small set of (unmeasurable) variables, called \emph{factor} variables.

These ideas have been exploited in so called \emph{dynamic factor} models, see e.g.  \cite{DEISTLER_1997}, \cite{deistler2007}, \cite{ForniHLR2000}, \cite{zorzi2015factor} and references therein.

Another possible avenue to control complexity is to build sparse dynamic models, e.g. exploring Granger's causality structure \citep{GRANGER_CAUSALITY} as done in \cite{CHIUSO_PILLONETTO_SPARSE_2012}.

In this paper we shall extend  and merge these ideas, building so called \emph{sparse plus low-rank (S+L)} models, \cite{BSL_CDC}; our aim is to model processes $y$  which Granger's causality structure 
is not necessarily  sparse, but it may become so, in an appropriate manner to be defined later on, after some latent variables (called \emph{factors} in analogy with factor models) are added. 

In this way, the relations among manifest variables and factor  variables will be described through a two-layer graph (or network) where the (few) nodes in the top (and hidden) layer denote the factor variables, whereas the ones in the bottom (and visible) layer denote the manifest variables.  The direct relation, in a sense to be precisely defined later on, between  two nodes is encoded by the presence of a connecting directed edge. If the relations among manifest variables are mostly encoded through the factor variables and the number of the latter is small (as compared to the number of manifest variables), 
this graph will be referred to as \emph{sparse plus low-rank (S+L) network}; in fact, as we shall see,  its structure translates into a S+L structure for the dynamic model of the manifest process. In particular, the rank of the  low-rank component will coincide with the   
 number of factor  variables whereas the sparse component depends on the number of edges  among the nodes representing the manifest variables. This modeling framework seems natural in many applications, such as biomedical sciences 
\citep{BRAIN_CDC}, econometrics  
\citep{LATENTG,Chandrasekaran_latentvariable,MULTIRES_CHOI_WILLSKY_2010}, and so on. Finally, this S+L structure have been considered also to perform robust principal component (PCA) analysis \citep{candes2011robust}.

In the first part of this paper we propose a new S+L network defined as a directed graph wherein edges encode conditional Granger causality dependences \citep{GRANGER_CAUSALITY,COND_GRANGER_CAUSALITY_2006} among variables. 
In the special case where there is no factor  variable our model coincides with the sparse model considered in \cite{CHIUSO_PILLONETTO_SPARSE_2012}. Moreover, the S+L model we propose is strictly connected to so called quasi-static factor models \citep{DEISTLER_1997,deistler2007,DEISTLER_2015}. Therefore, our S+L model can be understood as blend of those two models. 

Within this framework, we shall formulate 
an identification problem to select, using a finite set of measured data, the most appropriate S+L model (and thus the corresponding network). 

A consolidated paradigm in system identification is the so-called prediction error method (PEM), see  \cite{LJUNG_SYS_ID_1999,SODERSTROM_STOICA_1988}. In the traditional setting, candidate models are described in fixed parametric model structures, e.g. ARMAX. However, there are two main difficulties in this setting. First a model selection problem (i.e. order selection), usually performed by AIC and BIC criteria \citep{AKAIKE_1974,SCHWARZ_1978}, needs to be solved. Second, the parameterization of the predictor is non-linear,     so that minimizing the squared prediction error leads to a non-convex optimization problem. Regularization has been recently  introduced in the PEM framework, see
\cite{PILLONETTO_2011_PREDICTION_ERROR,PILLONETTO_DENICOLAO2010,EST_TF_REVISITED_2012,CHIUSO_PILLONETTO_SPARSE_2012,KERNEL_METHODS_2014,ChiusoARC2016}, as an alternative approach to control complexity of the estimated models. With this latter setting we search the candidate model, described via the predictor impulse responses, in an 
infinite dimensional nonparametric model class; the inverse (ill-posed) problem of determining a specific model using a finite set of measured data can be made into a well posed one  using a penalty term,  whose duty is to favor models with specific features. In the Bayesian view,  this is equivalent to the introduction of an {\em a priori} probability (prior)  on the unknown model. In  the nonparametric Gaussian regression approach proposed in \cite{PILLONETTO_2011_PREDICTION_ERROR}, this prior distribution is completely characterized by the covariance function, known also as kernel in the machine learning literature.   The kernel encodes the {\em a priori} knowledge   about the predictor impulse responses. In our case, the {\em a priori} knowledge is that the
predictor impulse responses must be Bounded Input Bounded Output (BIBO) stable and respect the S+L structure. Starting from these {\em a priori} assumptions, we derive the corresponding kernel by using the maximum entropy principle. In particular, we consider two possible alternative formulations to endow this {\em a priori} properties in the kernel, and thus  two different identification algorithms.

These kernels are characterized by the decay rate of the predictor impulse responses, by the number of conditional Granger causality relations among the manifest variables and by the number of factor  variables. This ensemble of features is not known and is  characterized by the so called hyperparameters vector. The latter is usually estimated by minimizing the negative log-marginal likelihood of the measured data \citep{RASMUSSEN_WILLIAMNS_2006}. In our case, the challenge is to perform the joint estimation of the hyperparameters tuning the sparse part and the low-rank one. Indeed, it should be observed that the S+L decomposition of a given model  might not be unique. On the other hand, once the hyperparameters vector is fixed the uniqueness of the estimated model will be guaranteed through regularization. To estimate the hyperparameters minimizing the negative log-marginal likelihood, we propose an algorithm imposing an 
``hyper regularizer'' on the low-rank hyperparameter to partially handle the non-uniqueness of the S+L decomposition.

Numerical experiments involving both S+L and generic models  show the effectiveness of our identification procedure both in terms of complexity, predictive capability and impulse response fit. 

The outline of the paper follows. In Section \ref{section_SL}, we introduce the S+L models. In Section \ref{section_pb_formulation}, we introduce the S+L identification problem. In Section \ref{section_kernel}, we derive the maximum entropy kernels inducing BIBO stability, sparsity and low-rank. Section \ref{section_SL_procedure} deals with the estimation of the hyperparameters vector.
In Section  \ref{section_simulation}, we provide some numerical examples to show the effectiveness of our method. Finally,  conclusions are drawn in Section
\ref{section_conclusions}. In order to streamline the presentation all proofs are deferred to the Appendix.

\subsection*{Notation}
Throughout the paper, we will use the following notation. $\Ns$ is the set of natural numbers. Given a finite set $\Ic$,
$| \Ic |$ denotes its cardinality.
$ \Es[\cdot]$ denotes the expectation, while $ \Es[\cdot|\cdot]$ denotes the conditional mean. Given three (possibly infinite dimensional) random vectors $\aaa$, $\bbb$ and $\ccc$ we say that $\aaa$ is conditionally independent of $\bbb$ given $\ccc$ if 
$$
 \Es[\aaa|\bbb,\ccc] =  \Es[\aaa|\ccc].
$$
Given $G\in\Rs^{n\times p}$, $[G]_{ij}$ denotes the entry of $G$ in position $(i,j)$. $\Mc^p$
denotes the vectors space of symmetric matrices of dimension $p$. $\Mc^p_+$ is the cone of symmetric positive definite matrices of dimension $p$, and 
$\overline{\Mc}^p_+$
denotes its closure. $\ell_2(\Ns)$ denotes the space of $\Rs$-valued infinite length sequences, which we think as infinite dimensional column vectors 
$g:=[g_1 \; g_2 \; \ldots \; g_j \; \ldots]^\top$,  $g_k\in\Rs$, $k\in\Ns$, such that
$\|g\|_2:=\sqrt{\sum_{k=1}^\infty |g_k|^2}<\infty$. 
$\ell_2^{p\times n}(\Ns)$ is the space of matrices of sequences in $\ell_2(\Ns)$
\al{\Phi=\left[\begin{array}{ccc} (\phi^{[11]})^\top & \ldots  & (\phi^{[1n]})^\top \\ \vdots  & \ddots & \vdots \\ (\phi^{[p1]})^\top & \ldots &  (\phi^{[pn]})^\top\end{array}\right]}
where $\phi^{[ij]}\in\ell_2 (\Ns)$, $i=1\ldots p$ and $j=1\ldots n$. 
$\ell_1(\Ns)$ denotes the space of $\Rs$-valued infinite length sequences $g$  
such that $\|g\|_1:=\sum_{k=1}^\infty |g_k|<\infty$. $\ell_1^{p\times n}(\Ns)$ is defined in similar way. 
$\Sc_2(\Ns)$ denotes the space of symmetric infinite dimensional matrices $K$ such that  
$\| K\|_2:=\sqrt{\sum_{i,j=1}^\infty | [K]_{ij}|^2}<\infty$. $\Sc_2^{p}(\Ns)$ is the space of symmetric infinite dimensional matrices 
\al{K=\left[\begin{array}{cccc} K^{[11]} & K^{[12]} & \ldots  & K^{[1p]} \\ 
K^{[12]} & K^{[22]} &  & K^{[2p]} \\ 
\vdots  & & \ddots & \vdots \\ K^{[1p]} & K^{[2p]} &\ldots & K^{[pp]}\end{array}\right]}
where $K^{[ij]} \in\Sc_2(\Ns)$, $i,j=1\ldots p$.
Given $\Phi\in\ell_2^{p \times n}(\Ns)$, $\Psi\in\ell_2^{m \times n}(\Ns)$ and $K\in\Sc_2^{n}(\Ns)$, the products $\Phi \Psi^\top$ and $\Phi K \Psi^\top$ are understood as $p \times m$ matrices whose entries are limits of infinite sequences \citep{INFINITE_MATRICES}. Given $g\in \ell_2^{n\times 1}(\Ns)$ and $K\in \Sc_2^{n}(\Ns)$,  $\|g\|^2_K:=g^\top K^{-1}g$. The definition is similar in the case that $g$ and $K$ have finite dimension. With some abuse of notation the symbol $z$ will denote both the complex variable as well as the shift operator 
$z^{-1} y(t) :=y(t-1)$.
Given a transfer matrix $L(z)$, $z\in\Cs$, of dimension $p\times p$, with some abuse of  terminology, we say that $L(z)$ has rank $n$, with $n\leq p$, if it admits the decomposition $L(z)=FH(z)$ where $F$ is a $p\times n$ matrix and $H(z)$ is a $n\times p$ transfer matrix. 
Given a stochastic process $y=\{y(t)\}_{t\in\Zs}$, its $i$-th component is denoted by $y_i = \{y_i(t)\}_{t\in \Zs}$.
With some abuse of notation, $y(t)$ will both denote a random vector and its sample value. 
From now on the time $t$ will denote \emph{present} and we shall talk about \emph{past} and \emph{future} with respect to time $t$. With this convention in mind,  \al{\yv^- =\left[\begin{array}{ccc}  y(t-1)^\top & y(t-2)^\top & \ldots  \end{array}\right]^\top}
denotes the (infinite length) past data vector of $y$ at time $t$. In similar way, $\yv^-_i$ denotes the past data vector of $y_i$ at time $t$.

\section{Sparse plus Low rank Models} \label{section_SL}
Consider a zero-mean stationary and Gaussian stochastic process $y$ taking values in $\Rs^p$.
Let $y$ be manifest, i.e. it can be measured, and described by the innovation  model
\al{  \label{OEmodel}y(t)=& G(z)y(t)+e(t)}
where $G(z)=\sum_{k=1}^\infty G_k z^{-k}$ and $(I-G(z))^{-1}$ are  BIBO stable transfer matrices of dimension $p\times p$ and  $e$ is  a zero mean white Gaussian noise (WGN) with  covariance matrix $\Sigma$.
The minimum  variance one-step ahead predictor of $y(t)$ based on the past data $\yv^-$, denoted by $\hat y(t|t-1)$ is given by 
\al{ \label{y_hat} \hat y(t|t-1)&= G(z)y(t)} so that $e$ is the one-step-ahead prediction error 
$$
e(t) = y(t) - G(z) y(t) = y(t) - \hat y(t|t-1).
$$
 Often $G(z)$ is approximated by a ``simple structure'' in order to simplify the analysis of the underlying system. One possible simple structure is $G(z)$ sparse, i.e. many of its entries are null transfer functions  
\citep{CHIUSO_PILLONETTO_SPARSE_2012}. Sparsity of the predictor transfer function $G(z)$ encodes the Granger's causality \citep{GRANGER_CAUSALITY} of $y$ and can be graphically represented with a   Bayesian network (see below).
However, in many interesting applications \citep{LATENTG,SOCIAL_NETWORK,SYSTEM_BIOLOGY,ForniHLR2000}   the components in $y$ are strongly related through common factors,  not accessible for observation. In such situations   approximating  $G(z)$ with a sparse matrix will likely yield to poor  results.

To overcome this limitation we consider a more general, but still ``simple'', structure by introducing an $n$ dimensional process $x$ ($n\ll p$), called   factor. We assume that $x$ is  zero-mean, jointly stationary and Gaussian with $y$. The underlying idea is that the (few) components of the factor $x$ are able to capture the common movement of $y$ so that, conditionally on $x$, the predictor of $y$ can now be well approximated by a sparse model. This can be formalized assuming that 
  $y$ can be modeled as follows:
\al{\label{S+Lmodel} y(t)=Fx(t)+S(z)y(t)+w(t)}
where $F\in\Rs^{n\times p}$ is the so-called \emph{factor loading} matrix and, in analogy with dynamic factor models, $Fx(t)$ is called the ``common component'' or ``latent variable'' of $y$  \citep{deistler2007}; $S(z)=\sum_{k=1}^\infty S_k z^{-k} $ is a BIBO stable and sparse transfer matrix; $w$ is a white noise process uncorrelated with the past histories of $y$ and $x$. As a result, the one-step ahead predictor of $y$ is given by
 \al{\label{S+Lpred}\hat y(t|t-1)&= F\hat x(t|t-1)+S(z)y(t)}
where $\hat x(t|t-1)$ is the minimum variance estimator of $x(t)$ based on $\yv^-$, that is 
\al{\label{xhat}\hat x(t|t-1): =  \Es[x(t)|  \yv^-] =H(z) y(t)} 
where $H(z)=\sum_{k=1}^\infty H_k z^{-k}$ is a  BIBO stable transfer matrix of dimension $n\times p$. We conclude 
that $\hat y(t|t-1)$ takes the form in (\ref{y_hat}) with $G(z)=S(z)+L(z)$ where $L(z):=FH(z)$. 
Hence, $G(z)$
has a sparse plus low-rank (S+L) structure.

\begin{figure}[htbp]
\begin{center}
\includegraphics[width=\columnwidth]{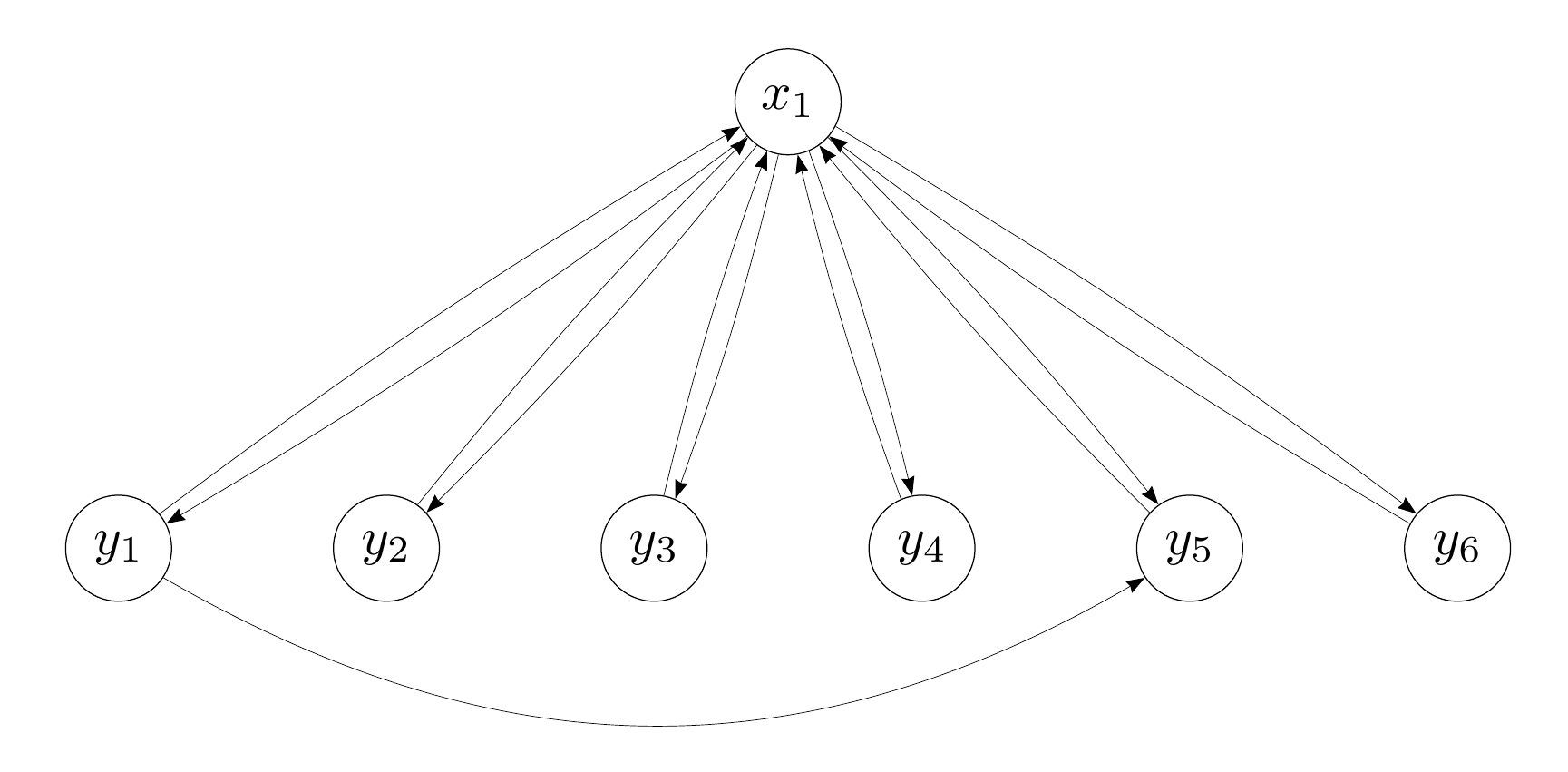}
\end{center}
\caption{Example of a network describing a S+L model with $p=6$ manifest variables and $n=1$ factor variable.}\label{graph_ex}
\end{figure}

It is possible to describe the structure of model  (\ref{S+Lmodel}) using a Bayesian network 
\citep{LAURITZEN_1996} 
having two layers. The nodes in the top layer correspond to the scalar processes $x_i$ $i=1\ldots n$, i.e. the components of $x$, whereas the ones in the bottom layer correspond to the scalar processes $y_i$ $i=1\ldots p$, i.e. the components of  $y$, see Figure \ref{graph_ex}.  
Then, the connections among the nodes obey the following rules:
\begin{itemize}
\item there is a  directed link from  node  $y_j$ to  node  $y_i$ if
$$
 \Es[y_i(t) |x(t),  \yv^-] \neq  \Es[y_i(t) |x(t), \yv_k^-, k=1,..,p, k\neq j]
$$
i.e. if $\yv_j^-$ is needed to predict $y_i(t)$  given $\yv_k^-$ with $k\neq j$ and $x(t)$. In this case, we shall say $y_j$ conditionally Granger causes $y_i$ 
\item there is a directed link from  node  $y_j$ to  node $x_i$ if
$$
\Es[x_i(t) | \yv^-] \neq \Es[x_i(t) | \yv_k^-, k=1,..,p, k\neq j]
$$
i.e. if  $\yv_j^-$ is needed to predict $x_i(t)$  given $\yv_k^-$ with $k\neq j$. In this case, we shall say $y_j$ conditionally Granger causes $x_i$
 \item there is a direct link from  node  $x_j$ to  node  $y_i$ if
 $$
 \Es[y_i(t) |\yv^-,x(t)] \neq  \Es[y_i(t) |\yv^-, x_k(t), k=1,..,n, k\neq j]
$$
i.e. if 
 $x_j(t)$ is needed to predict $y_i(t)$  given $\yv^-$ and $x_k(t)$ with $k\neq j$. In this case, we shall say $x_j$ conditionally Granger causes $y_i$.   
\end{itemize} 
Therefore, the S+L model in (\ref{S+Lmodel}) represents a network wherein manifest variables Granger cause each other mostly through  few factor variables. In Figure \ref{graph_ex}
we provide an example of a network describing a S+L model with $p=6$ manifest variables and $n=1$ factor variable. In particular, $y_5$ is conditionally Granger caused by $x$ and $y_1$, $y_i$ with $i\neq 5$ is conditionally Granger caused only by $x$, and $x$ is conditionally Granger caused by $y$. Therefore, the corresponding S+L model in (\ref{S+Lmodel}) has $S(z)$ with only one nonnull entry in position (5,1)
and $L(z)$ has rank equal to one.   
 
The decomposition of a transfer matrix $G(z)$ into sparse plus low-rank, i.e. $G(z)=S(z)+L(z)$ may not be unique.
As noticed in \citep{Chandrasekaran_latentvariable}, this degeneracy may occur when $L(z)$ is sparse, i.e. the factor  variables are not sufficiently ``diffuse'' across the manifest variables, or the degree of sparsity of $S(z)$ is low, i.e. there are manifest variables conditionally Granger caused  by too many manifest variables. On the other hand, it is possible to derive conditions under which such a decomposition is locally identifiable. We do not tackle this nonidentifiability issue in this paper, although it  is important.  Indeed, our aim is to find one
S+L decomposition (see the next Section) which is not necessarily unique.

\subsection{Connections: feedback models and quasi-static factor models}
It is also interesting to consider a ``feedback'' representation of the joint process $v:=[y^\top \; x^\top]^\top$ and see how this connects with S+L model \eqref{S+Lmodel} and other models already considered in the literature, such as quasi-static factor models \citep{deistler2007}.

Recall that  \citep{Gevers,Caines}, given the jointly stationary process  $v$, there is an essentially unique 
feedback model of the form 
\al{\label{dynamicfeedbackmodel}
\begin{array}{rcl}
y(t)&=&F(z)x(t)+d(t) \\
x(t)&=&C(z) y(t)+r(t)
\end{array}
}
where the driving noises  $r(t)$ and $d(t)$ are uncorrelated and jointly stationary processes (possibly non-white) and the feedback interconnection \eqref{dynamicfeedbackmodel} is internally stable, namely 
$$\left[\begin{array}{cc}
(I-F (z)C(z))^{-1} & (I-F(z) C(z))^{-1} F(z) \\ 
(I-C(z)F(z))^{-1} C(z) & (I-C(z)F(z))^{-1} 
\end{array}\right]
$$
 is analytic inside the closed unit disc and $C(\infty)=0$.  In particular, given an analytic and minimum phase spectral factor
$$
W(z)=\left[\begin{array}{cc}
W_{11}(z)&W_{12}(z)\\ 
W_{21}(z) & W_{22}(z)
\end{array}\right]
$$ of the joint spectrum $\Phi_v(z) = W(z) Q W^\top(z^{-1})$ of  $v$, 
lower block triangular and normalized at infinity (i.e. $W_{11}(\infty)=I$, $W_{22}(\infty)=I$, $W_{21}(\infty)=0$), the transfer function $F(z)$ is given by (see e.g. eq. (3.4) in \cite{Gevers}) $F(z) = W_{12}(z) W_{22}^{-1}(z)$.

The predictor for $y(t)$ given the past $\yv^-$ can be written, using \eqref{dynamicfeedbackmodel},  as 
$$
\hat y(t|t-1) = F (z)\hat x(t|t-1)  +\underbrace{ \hat\Es[d(t)| \yv^-]}_{\hat d(t|t-1)}.
$$
Hence, model  \eqref{S+Lmodel} is recaptured when 
\begin{enumerate}
\item[(a)] $
F(z) := W_{12}(z) W^{-1}_{22}(z) = F$ (i.e. $F(z)$ is  constant)  
\item[(b)] $\hat d(t|t-1) := S(z) y(t)$ is sparse.
\end{enumerate}
Note that condition (a) above is equivalent to $W_{12}(z) = F W_{22}(z)$, while complete freedom is left to the other entries of $W(z)$. 

To summarize, our model \eqref{S+Lmodel} is equivalent to assuming that there is a process $x$, jointly stationary with $y$, so that the pair $(y,x)$ admits the internally stable feedback representation:
\al{\label{feedbackmodel}
\begin{array}{rcl}
y(t)&=&Fx(t)+d(t) \\
x(t)&=&C(z) y(t)+r(t)
\end{array}
}
with 
uncorrelated (and possibly non-white) noises  $r(t)$ and $d(t)$.

 It is worth noting that when $C(z)=0$, model (\ref{feedbackmodel}) can be understood as a quasi-static factor model \citep{deistler2007} 
where $F$ is the factor loading matrix, $x$ is the $n$-dimensional factor and $d$ is colored noise uncorrelated with $x$, such that its minimum variance one-step ahead predictor based on $\yv^-$ is sparse: $\hat d(t|t-1)=S(z)y(t)$.

In general, however, the feedback model \eqref{feedbackmodel} is a generalization of quasi-static feedback models as $d(t)$ might be correlated with $x(t)$ (which happens when $C(z)\neq 0$ in \eqref{feedbackmodel}). 
The predictor model \eqref{S+Lpred} can be obtained from \eqref{feedbackmodel} simply  computing the predictor $\hat x(t|t-1)$, see \eqref{xhat}:
$$
\begin{array}{rcl}
\hat x(t|t-1) &=& C(z) y(t) + \hat r(t|t-1) \\ &=& C(z) y(t) + G_r(z) y(t) \\
&=& \underbrace{H(z)}_{C(z) + G_r(z) } y(t)
\end{array}
$$
where $G_r(z)y(t)= \hat r(t|t-1) $.

 \section{Problem Formulation}\label{section_pb_formulation}
Assume measured data $\{y(t)\}_{t=1,..,N}$ are available from the manifest process $y$ generated by (\ref{S+Lmodel}). 
The factor process $x$ cannot be measured nor its dimension $n$ is known. In this section, we address the problem of  estimating $S(z)$ and $L(z)$ from $\{y(t)\}_{t=1,..,N}$.

The transfer matrix $S(z)$ is parameterized in terms of its impulse response coefficients $S_k$. In particular, defining 
 $s^{[ij]}\in \ell_1(\Ns)$ to be the impulse response from input $j$ to output $i$ of the transfer matrix $S(z)$,  we have:  
  \al{  \label{def_s_ij} s^{[ij]}:= \left[
                  \begin{array}{ccccc}
                    [S_1]_{ij} & [S_2]_{ij} & \ldots & [S_k]_{ij}&\ldots \\
                  \end{array}
                \right]^\top . }
The coefficient vector $\theta_s^\top\in\ell_1^{1\times p^2}(\Ns)$ is defined  as follows:
\al{ \label{def_theta_s} \theta_s^\top =& \left[
                   \begin{array}{ccc|c}
                    (s^{[11]})^\top & \ldots  & (s^{[1p]})^\top & \ldots\\
                   \end{array}
                 \right.\nn\\
                 & \hspace{0.2cm} \left.
                   \begin{array}{c|ccc}
                      \ldots & (s^{[p1]})^\top & \ldots &(s^{[pp]})^\top \\
                   \end{array}
                 \right].}

Similarly, the impulse response coefficients $L_k$ parameterizing the transfer matrix $L(z)=\sum_{k=1}^\infty L_k z^{-k}$ are stacked in $\theta_l^\top \in \ell_1^{1\times p^2}(\Ns)$ as done above for $S(z)$. We first introduce   $l^{[ij]}\in\ell_1(\Ns)$ as 
 \al{ \label{def_l_ij} l^{[ij]}:= \left[
                  \begin{array}{ccccc}
                    [L_1]_{ij} & [L_2]_{ij} & \ldots  & [L_k]_{ij} \ldots \\
                  \end{array}
                \right]^\top. }
and define 
  \al{ \label{def_theta_l} \theta_l^\top =& \left[
                   \begin{array}{ccc|c}
                    (l^{[11]})^\top & \ldots  & (l^{[1p]})^\top & \ldots  \\
                   \end{array}
                 \right.\nn\\
                 & \hspace{0.2cm}
                 \left. \begin{array}{c|ccc}
                    \ldots & (l^{[p1]})^\top & \ldots &(l^{[pp]})^\top \\
                   \end{array}
                 \right].
                 }
 
The measured data $y(1)\ldots y(N)$ are stacked in the vector $\yv^+$ as follows
  \al{ \label{y_vettore}\yv^+ =&\left[
          \begin{array}{ccc|c}
            y_1(t)^\top &  \ldots & y_1(t+N-1)^\top  & \ldots  \\
          \end{array}
        \right.\nn\\
        & \hspace{0.2cm}\left.
          \begin{array}{c|ccc}
            \ldots  & y_p(t) ^\top&  \ldots & y_p(t+N-1)^\top \\
          \end{array}
        \right]^\top.}
         The vector $\ev^+$ is defined analogously.
Let us also introduce the {\em Toeplitz} matrix  $\phi_j \in\ell_2^{N\times 1}(\Ns)$:
 \al{ [\phi_j]_{kh}:=y_j(t+k-h-1)}
with $k=1\ldots N$ and $h\in\Ns$.

Then, we define the  regression matrix   $\Phi\in\ell_2^{pN \times p^2} (\Ns)$ as:
 \al{ \Phi=I_p \otimes \left[
                       \begin{array}{ccc}
                         \phi_1 & \ldots & \phi_p \\
                       \end{array}
                     \right]}   
so that, from (\ref{S+Lmodel}) the vector $\yv^+$ containing the measured data satisfy the linear model\footnote{ The matrix $\Phi$ contains in principle data in the remote past which are not available.  
In practice model \eqref{LM} needs to be approximated truncating $\Phi$ (and thus $\theta_l$ and $\theta_s$). This corresponds to assuming zero initial conditions, and it is a reasonable approximation given the  decay, as a function of $k$, of the coefficients $L_k$ and $S_k$ (which is a necessary condition for BIBO stability). We shall not enter into these details in the paper. See also Remark \ref{rem_regressor}.}  \al{\label{LM} \yv^+=\underbrace{\Phi(\theta_l+\theta_s)}_{\hat \yv^+}+\ev^+ }
where $\hat\yv^+:= \Phi(\theta_s+\theta_l)$ is the one-step ahead predictor of $\yv^+$.

Therefore, our S+L identification problem can be formulated in terms of PEM as follows.
 \pb \label{problem}Find $\theta_s^\top,\theta_l^\top\in\ell_1^{1\times p^2}(\Ns)$ corresponding to a S+L model  minimizing the prediction error norm $\| \yv^+-\Phi(\theta_s+\theta_l)\|^2_{\Sigma^{-1}\otimes I_N}$. \epb Following the nonparametric Gaussian regression approach in \cite{PILLONETTO_2011_PREDICTION_ERROR},  
we model $\theta_s$ and $\theta_l$ as two zero-mean processes
with kernels $K_S\in\Sc_2^{p^2}(\Ns)$ and $K_L\in\Sc_2^{p^2}(\Ns)$, respectively. These kernels may depend upon some tuning parameters, usually called hyperparameters and denoted with  $\bar \xi$ hereafter. As illustrated in Section \ref{section_kernel},  according to the maximum entropy principle $\theta_s$ and $\theta_l$ will be modeled as  Gaussian and independent.
In the following $\Hc_S$ and $\Hc_L$ denote the reproducing Hilbert spaces \citep{ARONSZAJN1950} of deterministic functions on $\Ns$, associated with $K_S$ and $K_L$, with norm denoted by $\|\cdot \|_{K_S^{-1}}$ and $\|\cdot \|_{K_L^{-1}}$, respectively. We assume that the past data $\yv^-$ neither affects the {\em a priori} probability on $\theta_s$ and $\theta_l$ nor carries information on $\bar \xi$ and $\Sigma$ \citep{POGGIO_GIROSI_1990}, that is 
\al{ \label{approx_pb}\mathbf{p} & (\yv^+,\theta_l,\theta_s,\yv^-|\bar\xi,\Sigma)\nn\\ 
& =\mathbf{p}(\yv^+|\theta_l,\theta_s,\yv^-,\bar \xi,\Sigma) \mathbf{p}(\theta_l,\theta_s|\yv^-,\bar \xi,\Sigma)\mathbf{p}(\yv^-|\bar \xi,\Sigma)\nn\\
& \approx
\mathbf{p}(\yv^+|\theta_l,\theta_s,\yv^-,\bar\xi,\Sigma) \mathbf{p}(\theta_l,\theta_s|\bar \xi,\Sigma)\mathbf{p}(\yv^-).}
Let 
\al{\hat \theta_s= \Es[\theta_s|\yv^+,\bar \xi,\Sigma],\;\; \hat \theta_l= \Es[\theta_l|\yv^+,\bar \xi,\Sigma]}  be, respectively, the  minimum variance estimators of $\theta_s$ and $\theta_l$  given $\yv^+$, $\bar \xi$ and $\Sigma$. In what follows, $\hat S(z)$ and $\hat L(z)$ denote the transfer matrices corresponding to $\hat \theta_s$ and $\hat \theta_l$, respectively.
The next Proposition shows that $\hat \theta_s$ and $\hat \theta_l$ are,  almost surely, solution to a {\em Tikhonov}-type variational problem and belong to the spaces $\Hc_S$ and $\Hc_L$, respectively. \pro   \label{prop_theta_SL}
Under the assumption that $y$ is a second order stationary  process and under approximation (\ref{approx_pb}),
almost surely we have
\al{ \label{Tickho_SL}(\hat \theta_s,\hat \theta_l)=\underset{\substack{\theta_s\in\Hc_{K_S}\\ \theta_l\in\Hc_{K_L}}}{\arg\min}\|\yv^+-\Phi &  (\theta_s+\theta_l)\|^2_{\Sigma^{-1}\otimes I_N}\nn\\
 & +\| \theta_s\|^2_{K_S^{-1}}+\| \theta_l\|^2_{K_L^{-1}}.}
Moreover, almost surely: 
\al{ \label{Bayes_SL}\hat \theta_s= K_S\Phi^\top c,\;\; \hat \theta_l=K_L \Phi^\top c   }
where \al{ c=(\Phi(K_S+K_L)\Phi^\top+\Sigma\otimes I_N)^{-1}\yv^+.} \epro

\rema \label{rem_regressor}The semi-infinite regression matrix $\Phi$ depends on both $\yv^+$ and $\yv^-$. Since $\yv^-$ is never completely known, a solution to handle the initial conditions consists of setting its unknown components to zero. In this way, the introduced error goes to zero as $N$ increases
 \citep[Section 3.2]{LJUNG_SYS_ID_1999}.  Alternatively it would be possible to incorporate initial conditions in the estimation problem, e.g. modeling also the free response of the system. This is however outside the scope of the paper and is only practically relevant when very slow dynamics are present.\erema 
The main task now is to design the kernels $K_S$ and $K_L$
in such a way that  $\hat S(z)$ and $\hat L(z)$ are almost surely BIBO stable while favouring $\hat S(z)$ to be sparse and $\hat L(z)$ to be of low-rank.

 \section{Maximum entropy priors} \label{section_kernel}
One way to derive a  probability law for the joint process $[\theta_s^\top \; \theta_l^\top]^\top$ under  desired constraints  rests on the maximum entropy principle. 
The most common justification of maximum entropy solutions relies on  ``information'' arguments  essentially stating that the maximum entropy distribution is the one  which entails the maximum ``uncertainty'' under the given constraints. There is another and very important motivation for adopting the maximum entropy solution:   \cite{SHORE_AXIOMATIC_1980}  have shown that maximum entropy is the unique correct method satisfying some minimal consistency axioms;  basically, these axioms state that the solution should be consistent when ``there are different ways of taking the same information into account''.

 We shall make the rather mild assumption that the process $[\theta_s^\top\; \theta_L^\top]^\top$ is zero-mean and absolutely continuous with respect to the Lebesgue measure. We will see that, under suitable constraints,  the optimal solution (i.e. maximizing the differential entropy) is a Gaussian 
process where $\theta_s$ and $\theta_l$ are independent. Then, we will also characterize the corresponding kernels. 
In what follows, we propose two different ways to enforce BIBO stability and low-rank on  $\hat L(z)$. This leads to two different types of kernel for $\theta_l$.

\subsection{First type of kernel}\label{first_type_kernel}
We start with the constraints on $\theta_s$ inducing BIBO stability and sparsity on $\hat S(z)$. To do so, we exploit the following proposition:

  \pro \label{prop_stability} Let $P\in\Sc_2(\Ns)$ be a strictly positive definite kernel (in the sense of Moore) such that $[P]_{tt}\leq \kappa t^\alpha e^{-\beta t}$, $t\in\Ns$, with $\kappa,\beta>0$ and $\alpha\in\Rs$.  Let also $\phi$ be a zero-mean process which satisfies the moment  constraint
    \al{\label{condps}\Es[\| \phi\|_{P^{-1}}^2]\leq c} 
  where $c\geq 0$.  Then, for any $\varepsilon>1$ there exists $\bar \kappa_\varepsilon>0$ such that   the covariance function (kernel) $K$ of  $\phi$ satisfies:
 \al{[K]_{tt}\leq \bar \kappa_\varepsilon t^{\alpha+\varepsilon} e^{-\beta t},\;\; t\in\Ns.} 
  \epro 
 
Thus,  we consider the constraint on $\theta_s$ \al{ \label{constraint_p_s}  \Es[\| s^{[ij]}\|^2_{P^{-1} }] \leq c_{ij},}
where $c_{ij} \geq 0$, $i,j=1\ldots p$, and $P\in\Sc_2(\Ns)$.
If $c_{ij}=0$, then $s^{[ij]}$ is the null sequence  in mean square and
 so is its posterior mean.

It is not difficult to see that the assumption on $P$ in Proposition \ref{prop_stability} is satisfied by the kernels usually employed in the identification of dynamical models (e.g. stable spline, tuned/correlated and so on, see \cite{KERNEL_METHODS_2014}). Thus, by (\ref{constraint_p_s}) the covariance of the $k$-th element of $s^{[ij]}$ decays exponentially. 
We conclude that the posterior mean of the transfer function in position $(i,j)$ of $S(z)$, under constraint (\ref{constraint_p_s}), is BIBO stable. Moreover, it is  null if and only if  $c_{ij}=0$. 

\begin{rem}
Clearly, if $P$ in \eqref{constraint_p_s} is chosen as the covariance matrix of DC (TC) prior, the resulting Maximum Entropy prior coincides with the DC (TC) prior. Some recent literature has discussed the maximum entropy properties of some kernels, such as the TC or DC kernels, see \cite{NicFer:1998:IFA_1779},  \cite{CarliCL2015} and \cite{ChenACCLP15}, as well as  extensions to more articulated kernels \cite{PrandoPCAuto2015}.   In these papers it was shown for TC and DC kernels that weaker  constraints with respect to \eqref{constraint_p_s} lead to the same solution.
We refer the reader to these papers form more details; here we stick to the simpler conditions   \eqref{constraint_p_s} in order to streamline the derivation as well as to allow alternative choices of the matrix $P$ in \eqref{constraint_p_s}.
\end{rem}

In view of Proposition \ref{prop_stability}, to guarantee BIBO stability on $\hat L(z)$ we also impose the constraint 
\al{  \label{constraint_p_l1}\sum_{i,j=1}^p\Es[\| l^{[ij]}\|^2_{ P^{-1} }]\leq r}
for some $r$ such that $0<r<\infty$.

We now switch our attention to the low rank property of the matrix $L(z)$. Let $A_l$ be the random semi-infinite matrix defined as \al{\label{def_Al} A_l=\left[\begin{array}{ccccc}L_1 & L_2 & \ldots  & L_k &\ldots \end{array}\right]}
and consider the constraint on $\theta_l$
\al{ \label{constraint_p_l2}\Es[A_l A_l^\top]\leq Q}
with $Q\in\overline{\Mc}^p_+$. If $Q$ has $p-n$ singular values equal to zero, so does the covariance $\Es[A_l A_l^\top]$; therefore    the posterior mean $\hat A_l$ of $A_l$ has rank less than or equal to $n$ thus  admitting the decomposition
 \al{ \hat A_l=\left[\begin{array}{ccccc} \hat F \hat H_1 & \hat F\hat H_2 & \ldots & \hat F\hat H_k &\ldots \end{array}\right],}
where $\hat F \in \Rs^{p\times n}$ and the $\hat H_k \in \Rs^{n\times p}$, $k\in\Ns$, as in Section \ref{section_SL}. 
We conclude that, under constraints (\ref{constraint_p_l1}), $\hat L(z)$ is BIBO stable and under constraint (\ref{constraint_p_l2}), its rank is less than or equal to $n$
if and only if $Q$ has rank equal to $n$.

In order to build the desired prior distribution we make use of Kolmogorov extension Theorem, see  \cite{Oksendal}, and work with finite vectors extracted from processes $\theta_l$ and $\theta_s$.

To do so, let us consider a finite index set $\Ic=\Ic_s \times \Ic_l$ in $\Ns \times \Ns$. Let $\check \theta_s$ and $\check \theta_l$ be the random vectors whose components are extracted, respectively, from the process $\theta_s$ and $\theta_l$ according to the index sets $\Ic_s$ and $\Ic_l$. We denote by  $\pp_{\Ic}(\check \theta_s,\check\theta_l)$ the joint probability density of $\check \theta_s$ and $\check\theta_l$.
 By Kolmogorov extension Theorem the joint process $[ \theta_s^\top \; \theta_l^\top]^\top$ can be characterized by specifying the joint probability density $\pp_{\Ic}$ for all finite sets $\Ic \subset \Ns\times \Ns$. Thus, the maximum entropy process $[ \theta_s^\top \; \theta_l^\top]^\top$ can be constructed building all marginals   $\pp_{\Ic}$ using the maximum entropy principle, which can thus be extended by Kolmogorov extension theorem.
 Such principle states that
among all the probability densities $\pp_{\Ic}$ satisfying the desired constraints, the optimal one should maximize the differential entropy \citep{COVER_THOMAS}
 \al{ \mathbf{H}(\pp_{\Ic})=-\Es[ \log(\pp_\Ic)].}
Constraints (\ref{constraint_p_s}), (\ref{constraint_p_l1}) and (\ref{constraint_p_l2}) boil down, respectively, to
\al{ & \label{constr_s}\Es[\| \check s^{[ij]}\|^2_{P^{-1}_{\Ic_s} }]\leq c_{ij}, \;\; i,j=1\ldots p\\
& \label{constr_l1}\sum_{i,j=1}^p \Es[\| \check l^{[ij]}\|^2_{P^{-1}_{\Ic_l} }]\leq r\\
& \label{constr_l2} \Es[\check A_l \check A_l^\top]\leq Q}
where $\check s^{[ij]}$ and $\check l^{[ij]}$ are the vectors extracted from $s^{[ij]}$ and $l^{[ij]}$ according to the index set $\Ic_s$ and $\Ic_l$, respectively.  
$ P_{\Ic_s}$ and $ P_{\Ic_l}$ are the kernel matrices whose entries are extracted from $ P$ according to $\Ic_s$ and $\Ic_l$, respectively. $\check A_l$ is the matrix whose blocks are extracted from $A_l$ according to $\Ic_l$. 
Therefore, we obtain the following maximum entropy problem 
 \al{ \label{ME_problem}\underset{\pp_{\Ic}\in\Pc}{\max}& \mathbf{H} (\pp_{\Ic})\nn\\
\hbox{s.t. } & \Es[\| \check s^{[ij]}\|^2_{ P^{-1}_{\Ic_s} }]\leq c_{ij}, \;\; i,j=1\ldots p \nn\\
& \sum_{i,j=1}^p \Es[\| \check l^{[ij]}\|^2_{P^{-1}_{\Ic_l} }]\leq r\nn\\
& \Es[\check A_l \check A_l^\top]\leq Q}
 where $\Pc$ is the class of probability densities in $\Rs^{|\Ic_s|}\times\Rs^{ |\Ic_l|}$
which are bounded and taking positive values.

\teo\label{teo_ME} Under the assumption that $c_{ij}>0$, $i,j=1\ldots p$, and $Q\in \Mc^p_+$, the unique optimal solution to the maximum entropy problem (\ref{ME_problem}) is 
such that 
$\check \theta_s$ and $\check \theta_l$ are independent, Gaussian with zero mean and kernel matrix, respectively,
\al{  \label{optimal_Ks}\check K_S&= \mathrm{diag}(\gamma_1\ldots \gamma_{p^2})\otimes P_{\Ic_s}\\
 \label{optimal_Kl} \check K_L&=(\lambda I_{p^2}\otimes P_{\Ic_l}^{-1} +\tilde \Lambda \otimes I_{p|\Ic_l |})^{-1}}
 where $\gamma_i>0$, $i=1\ldots p^2$, $\tilde \Lambda\in\overline{\Mc}_+^p$ and $\lambda\geq 0$.
  \eteo

In what follows we assume that constraints (\ref{constr_l1}) and (\ref{constr_l2}) are totally binding in problem (\ref{ME_problem}). Since $\lambda$ and $\tilde \Lambda$
are the Lagrange multipliers associated to those constraints, it is not difficult to see that in this situation $\lambda>0$ and $\tilde \Lambda\in\Mc^p_+$. Moreover, we define $\Lambda=\tilde \Lambda^{-1}$. As noticed before, we are interested in the limiting cases where $c_{ij}$ might be equal to zero  and $Q$ might be of  low-rank in order to have sparse and low-rank estimators. 
To include these scenarios, we consider the limits as $c_{ij} \rightarrow 0 $ and $Q$ tends to a low rank matrix and extend the maximum entropy solution by continuity. 

\pro  \label{extended_ME} Let $\mathbf{C}=\{(i,j) \hbox{ s.t. } c_{ij}=0\}$ and $\mathbf{Q}=\{v \hbox{ s.t. } Qv=0\}$.  Then, the maximum entropy solution extended by continuity is the probability density such that      
 $\check \theta_s$ and $\check \theta_l$ are independent, Gaussian, zero-mean, and with kernel matrices
\al{ \check K_S&=  \mathrm{diag}(\gamma_1\ldots \gamma_{p^2}) \otimes P_{\Ic_s}\\
\check K_L&=\lambda^{-1} I_{p^2}\otimes P_{\Ic_l}
-\lambda^{-2} I_{p^2}\otimes P_{\Ic_l} \nn\\ & \hspace{0.5cm} \times (\lambda^{-1} I_{p^2}\otimes P_{\Ic_l} +\Lambda \otimes I_{p|\Ic_l |} )^{-1} I_{p^2}\otimes P_{\Ic_l} }
 where $\gamma_{(i-1)p+j}=0$ if and only if $(i,j)\in \mathbf{C}$ and $\Lambda v=0$ if and only if $v\in\mathbf{Q}$.   
\epro 

Finally, in view of Kolmogorov extension Theorem, from the probability density of $[\check \theta_s^\top \; \check\theta_l ^\top]^\top$ we can characterize the probability law of $[\theta_s^\top \; \theta_l ^\top]^\top$ maximizing the differential entropy.

\corr Consider the joint process $[\theta_s^\top \;\theta_l^\top ]^\top$ where $\theta_s$ and $\theta_l$ are Gaussian independent processes 
with kernels, respectively,
\al{ \label{ME_kernel1} K_S&= \mathrm{diag}(\gamma_1\ldots \gamma_{p^2}) \otimes P\nn\\
K_L & =\lambda^{-1} I_{p^2}\otimes P
-\lambda^{-2} I_{p^2}\otimes P \nn\\ & \hspace{0.5cm} \times (\lambda^{-1} I_{p^2}\otimes P +\Lambda \otimes I_{\infty} )^{-1} I_{p^2}\otimes P}
where $I_\infty$ such that $[I_\infty]_{tt}=1$, $t\in\Ns$, and zero otherwise.   
For all finite sets $\Ic\subset \Ns\times \Ns$, its joint probability density is the extended solution to the maximum entropy problem (\ref{ME_problem}).\ecorr

It is worth noting that the maximum entropy kernels are characterized by the hyperparameters $\gamma_i$, $i=1\ldots p^2$,
which control sparsity on $\hat S(z)$, $\Lambda$ tuning the rank (and column space) of $\hat L(z)$, while $\beta$ (see Proposition \ref{prop_stability}) controls the decay rate (as a function of the time index $k\in \Ns$) of the  estimators  $\hat L_k$ and $\hat S_k$ and thus BIBO stability of $\hat S(z)$ and $\hat L(z)$, see \cite{PILLONETTO_DENICOLAO2010}. 
Finally, $\lambda$ represents a trade-off between BIBO stability and low-rank on $\hat L(z)$.
The structure of $\Lambda$ is very general. We suggest that $\Lambda$ can be reparameterized by introducing few hyperparameters to reduce its {\em degrees of freedom}. In Section \ref{section_SL_procedure} we will propose one possible reparameterization.

\rema \label{rem_wipf} It  is possible to derive the same structure for $K_L$ by adopting the regularization point of view.
In \cite{WIPF_2012}, it has been shown that the
penalty term  $\log \det (R)$, with $R\in\Mc_+^p$, induces low-rank on $R$. Moreover, the term $\log \det (R)$ admits the variational upper bound
\al{ \log \det( R)\leq \tr(\Lambda^{-1} R)-\log \det(\Lambda^{-1} )-p }
where $\Lambda\in\Mc_+^p$ and equality holds if and only if $R=\Lambda$. 
Consider the random vector  $\check \theta_l$ extracted from $\theta_l$ according to $\Ic_l$.
 Thus, we can induce low-rank on $R=\check A_l \check A_l^\top$ by considering the penalty
\al{ \label{quadratic_penalty_wipf}  \tr(\Lambda^{-1} \check A_l \check A_l^\top)-\log  \det(\Lambda)-p }
where $\Lambda$ represents a rough estimate of $\check A_l \check A_l^\top$.
The unique term depending on $\check\theta_l$
in (\ref{quadratic_penalty_wipf}) is $\tr(\Lambda^{-1} \check A_l  \check A_l^\top)=\| \check \theta_l\|_{\Lambda^{-1}\otimes I_{ p |\Ic_l|}}^2$ which is one part of the norm $\|\check \theta_l \|_{\check K_L^{-1}}^2$ with 
$\check K_L$ kernel matrix defined in (\ref{optimal_Kl}). 
Thus,
this penalty terms induces low-rank on $\check A_l$.
\erema

\subsection{Second type of kernel} \label{second_type_kernel}

The BIBO stability and low-rank constraints on $\hat L(z)$ can be imposed by using only one constraint. Consider the random semi-infinite matrix $A_l$ defined in (\ref{def_Al}) and the  constraint \al{\label{constraint_p_l3} \Es[A_l (P^{-1}\otimes I_p)A_l^\top]\leq Q}
with $Q\in\overline{\Mc}^p_+$. Similarly to the previous case, if the null space of $Q$  has dimension $m-n$, then the posterior mean $\hat A_l$ of $A_l$ has rank less than or equal to $n$; therefore also  $\hat L(z)$ has rank less than or equal to $n$. This statement is formalized in the following proposition:
 \pro \label{prop_stability2} Assume that $P\in\Sc_2(\Ns)$ is strictly positive definite and such that $[P]_{tt}\leq \kappa t^{\alpha} e^{-\beta t}$, $t\in\Ns$, with $\kappa,\beta>0$ and $\alpha\in\Rs$. Then,  under constraint (\ref{constraint_p_l3}), for any $\varepsilon> 1$ there exists $\bar \kappa_\varepsilon >0$ such that
 $\Es[| [L_t]_{ij} |^2]\leq \bar \kappa_\varepsilon t^{\alpha+\varepsilon} e^{-\beta t}$, $i,j=1\ldots p$, $t\in\Ns$. In addition, if the vector $v$ belongs to the nullspace of $Q$, then
 $v^\top A_l$ is zero in mean square;  therefore $A_l$ has low rank (in mean square) and its null space contains the one of $Q$.\epro

Similarly to the previous case, the joint process $[\theta_s^\top \; \theta_l^\top]^\top$ is characterized through $\pp_{\Ic}$ using the maximum entropy 
principle. In particular, constraint (\ref{constraint_p_l3}) becomes \al{ \Es[\check A_l (P_{\Ic_l}^{-1}\otimes I_p)\check A_l^\top]\leq Q. }
The corresponding maximum entropy problem is 
\al{ \label{ME_problem2}\underset{\pp_{\Ic}\in\Pc}{\max}& \mathbf{H} (\pp_{\Ic})\nn\\
\hbox{s.t. } & \Es[\| \check s^{[ij]}\|^2_{P^{-1}_{\Ic_s} }]\leq c_{ij}, \;\; i,j=1\ldots p \nn\\
& \Es[\check A_l (P_{\Ic_l}^{-1}\otimes I_p)\check A_l^\top]\leq Q}

\teo\label{teo_ME2} Under the assumption that $c_{ij}>0$, $i,j=1\ldots p$, and $Q\in {\Mc}_+^p$, the unique optimal solution to the maximum entropy problem (\ref{ME_problem2}) is 
such that 
$\check \theta_s$ and $\check \theta_l$ are independent, Gaussian with zero mean and kernel matrix, respectively,
\al{  \label{optimal_Ks2}\check K_S&= \mathrm{diag}(\gamma_1\ldots \gamma_{p^2})\otimes P_{\Ic_s}\\
 \label{optimal_Kl2} \check K_L&=\Lambda\otimes I_{p}\otimes P_{\Ic_l}}
 where $\gamma_i>0$, $i=1\ldots p^2$,
and  $\Lambda\in\Mc_+^p$. 
  \eteo

Similarly to the previous case, we extend the maximum entropy solution by continuity to the case of interest.

\pro \label{extended_ME2} Let $\mathbf{C}=\{(i,j) \hbox{ s.t. } c_{ij}=0\}$ and $\mathbf{Q}=\{v \hbox{ s.t. } Qv=0\}$.  Then, the maximum entropy solution extended by continuity is the probability density such that      
 $\check \theta_s$ and $\check \theta_l$ are independent, Gaussian, zero-mean, and with kernel matrices
(\ref{optimal_Ks2}) and (\ref{optimal_Kl2}) , respectively, where $\gamma_{(i-1)p+j}=0$ if and only if $(i,j)\in \mathbf{C}$ and $\Lambda v=0$ if and only if $v\in\mathbf{Q}$.  
\epro

Finally, from the probability density of $[\check \theta_s^\top\; \check \theta_l^\top]^\top$ we characterize the probability law of the joint process $[\theta_s^\top\; \theta_l^\top]^\top$ using the Kolmogorov extension Theorem.

\corr Consider the joint process $[\theta_s^\top \;\theta_l^\top ]^\top$ where $\theta_s$ and $\theta_l$ are Gaussian independent processes 
with kernel, respectively,
\al{ \label{ME_kernel2}K_S&=\mathrm{diag}(\gamma_1\ldots \gamma_{p^2})  \otimes P\nn\\
K_L &=  \Lambda  \otimes I_{p}\otimes P.}
For all finite sets $\Ic\subset \Ns\times \Ns$, its joint probability density is the extended solution to the maximum entropy problem (\ref{ME_problem2}). \ecorr

In respect to the first type of kernel, we do not need the extra hyperparameter $\lambda$. Furthermore, in order to compute $K_L$ in (\ref{ME_kernel2}) we do not need to invert an infinite dimensional matrix\footnote{Since the initial conditions are set equal to zero, in practice this corresponds to invert a matrix whose dimension is proportional to size of the data.} as in (\ref{ME_kernel1}). We conclude that the computation of $K_L$ in the second type of kernel is more efficient.

\rema The derivation of the  maximum entropy kernels given above requires that a matrix   $P$, satisfying the assumptions of Proposition \ref{prop_stability} and Proposition \ref{prop_stability2},   be fixed.
In what follows, we consider the filtered kernel proposed in \cite{PILLONETTO_2011_PREDICTION_ERROR}
\al{ \label{def_K_tilde} P= F P_{SS} F^\top} where 
$P_{SS}\in\Sc_2(\Ns)$ is the stable spline (SS) kernel
and $F$ is a Toeplitz matrix which can be used to shape the kernel emphasising certain frequencies. For instance, when
the focus is prediction, the predictor impulse  might exhibit an oscillatory behavior  due to ``high'' frequency zeros of the noise spectrum. In this paper,    following \cite{PILLONETTO_2011_PREDICTION_ERROR},  we have built $F$ from the 
 impulse response of a second order oscillatory  system; the poles of this system are estimated as hyperparameters.
  Note that, instead of SS one might choose other types of kernel, such as diagonal, tuned/correlated, diagonal/correlated and so on.
 \erema

 \section{Hyperparameters Estimation}\label{section_SL_procedure}

 In order  to compute $\hat \theta_s$ and $\hat \theta_l$, estimates of the noise variance  $\Sigma$ and of the hyperparameters vector $\bar \xi$ are needed. Let $\bar \xi:=\{ \tau,\xi\}$ where $\tau$ denotes the hyperparameters of $P$ defined in (\ref{def_K_tilde}), $\xi:=\{ \gamma_1\ldots \gamma_{p^2},\Lambda,\lambda\}$ for the first type of kernel, and  $\xi:=\{ \gamma_1\ldots \gamma_{p^2}, \Lambda\}$ for the second one.
In this paper $\Sigma$ is estimated  using a low-bias ARX-model as suggested in \citep{GOODWIN_1992}.
To estimate the hyperparameters of $P$ we consider the unstructured model for $y$
\al{\label{OEmodel} y(t)= G(t)y(t)+e(t)}
 where $G(z)=\sum_{k=1}^\infty G_k z^{-k}$ is BIBO stable and $e$ WGN with covariance matrix $\Sigma$.  
Equation \eqref{OEmodel} can be written as a linear regression model  \al{\label{oss_oe}\yv^+ =\Phi \theta +\ev^+} where $\theta^\top \in\ell ^{1\times p^2}(\Ns)$ contains the coefficients of $G(z)$. Following the Gaussian regression approach, $\theta$  in  \eqref{oss_oe} is modelled as a zero-mean Gaussian process with kernel $K=I_{p^2}\otimes P$. Doing so,   the hyperparameters of $P$
can be estimated by minimizing the negative log-marginal likelihood of $\yv^+$ computed using model  (\ref{oss_oe}), see \cite{PILLONETTO_2011_PREDICTION_ERROR}.

Then the hyperparameters vector  $\xi$ can be estimated minimizing
the negative log-marginal likelihood $\ell$ of $\yv^+$ under model  (\ref{LM}) with $P$ fixed as above. Under the assumptions of Proposition  \ref{prop_theta_SL}, we have 
 \al{\label{loglik}\begin{array}{rcl}
 \ell (\yv^+,\xi)&=&\frac{1}{2}\log \det V+\frac{1}{2}(\yv^+) ^\top V^{-1}\yv^++\hbox{ const. term}\\
 \mbox{} & & \\
 V&=& \Phi (K_S+K_L)\Phi^\top+\Sigma \otimes I_{N}
 \end{array}}
where $K_L$ and $K_S$ are defined in \eqref{ME_kernel1} or in \eqref{ME_kernel2}.
Notice that, the minimization of (\ref{loglik}) with respect to $\xi$ is a nonconvex constrained optimization problem. Accordingly, only
local minima can be computed. However, the real challenge is to perform the joint estimation of the $\gamma_i$'s and $\Lambda$ 
because the sparse and the low-rank part might be not identifiable from the measured data.
As a partial remedy it is useful  to reduce the \emph{degrees of freedom} of  $\Lambda$ (see also \cite{PrandoPCAuto2015}), constraining its structure as follows:
\al{ \label{Lambda_struct} \Lambda =\alpha (I-UU^\top)+U  \mathrm{diag} (\beta_1\ldots \beta_n) U^\top .}
The matrix  $U \in\Rs^{p\times n} $ is built, in view also of Remark \ref{rem_wipf}, using  the first $n$ singular vectors of an estimate $\hat A_l \hat A_l^\top$ of $A_lA_l^\top$, that is $U=[u_1\, \ldots\, u_n\,]^\top$ where $u_1\ldots u_n$ are the first $n$ singular vectors of  $\hat A_l \hat A_l^\top$. In this way, the 
constraints in $\Lambda$ are decoupled along the ``most relevant'' $n$ singular vectors of $\hat A_l\hat A_l^\top$ and their orthogonal complement. 
 This has the effect of steering the factor loading matrix $F$ toward the columns space of $U$.
Regarding the hyperparameters $\gamma_i$'s, it has been shown in the literature  \citep{McKayARD,Tipping01,
 ARAVKIM_CONVEX_NONCONVEX_2014} that the minimization of (\ref{loglik}) automatically leads to  sparsity in the $\hat\gamma_i$'s and therefore in $\hat S(z)$. Therefore, the minimization  of (\ref{loglik}) is performed with respect to $\tilde \xi:=\{\gamma_1\ldots \gamma_{p^2},\alpha,\beta_1\ldots \beta_n,\lambda\}$ for the first type of kernel and with respect to $\tilde \xi:=\{\gamma_1\ldots \gamma_{p^2},\alpha,\beta_1\ldots \beta_n\}$ for the second one while $n$ and $U$ are fixed.  
 The complete procedure 
 is outlined in Algorithm \ref{algo:hyper}. $n^{(k)}$, $U^{(k)}$, $\tilde \xi^{(k)}$, $\hat A_l^{(k)}$
and $\hat L^{(k)}(z)$ denote, respectively, $n$, $U$, $\tilde \xi$, $\hat A_l$
and $\hat L(z)$ at the $k$-th iteration of the algorithm. The marginal likelihood function \eqref{loglik}, when $\Lambda$ is constrained to be of the form \eqref{Lambda_struct}, is a function of $U$. This dependence is made explicit in the notation 
 \al{\label{loglikU}  \ell (\yv^+,\tilde \xi,U )&=\frac{1}{2}\log \det V+\frac{1}{2}(\yv^+) ^\top V^{-1}\yv^++\hbox{ const. term} \nn\\
V & =  \Phi (K_S+K_L)\Phi^\top+\Sigma \otimes I_{N} \nn\\
\Lambda & =  \alpha (I-UU^\top)+U  \mathrm{diag} (\beta_1\ldots \beta_n) U^\top\nn}
where $K_S$ and $K_L$ are defined in \eqref{ME_kernel1} for the first type of kernel and in \eqref{ME_kernel2} for the second one.
\begin{algorithm}
\caption{Computation of $n$, $U$ and $\tilde \xi$}
\label{algo:hyper}
\begin{algorithmic}[1] \small
\STATE $k=0$
\STATE $n^{(0)}\leftarrow 0$
\STATE $U_{OPT}^{(0)} \leftarrow $ empty matrix 
\STATE $\tilde \xi^{(0)} \leftarrow \underset{\tilde \xi}{\mathrm{argmin}}\,\ell(\yv^+,\tilde \xi,U_{OPT}^{(0)})$ 
\STATE $\tilde \xi_{OPT}^{(0)} \leftarrow \tilde \xi^{(0)}$
\REPEAT 
\STATE $k \leftarrow k+1$
\STATE $n^{(k)}\leftarrow n^{(k-1)}+1$
\IF{$n^{(k)}=1$}
\STATE $\hat A_l^{(k)}\leftarrow [\, \hat G_1\;\;  \hat G_2\;\; \ldots \,] $ where  $ \hat G_1, \hat G_2,\ldots$ are the
\STATEx \hspace{0.5cm}   coefficients of $\hat G(z)$ estimated from (\ref{OEmodel}) with 
\STATEx \hspace{0.5cm}  $K=I_{p^2}\otimes P$ 
\ELSE
\STATE $\hat A_l^{(k)}\leftarrow [\, \hat L_1^{(k)}\;\;  \hat L_2^{(k)}\;\; \ldots \,] $ where  $ \hat L_1^{(k)}, \hat L_2^{(k)},\ldots$ are 
\STATEx \hspace{0.5cm} the coefficients of $\hat L^{(k)}(z)$ estimated from (\ref{Tickho_SL}) with 
\STATEx \hspace{0.5cm}  $K_S$ and $K_L$ having hyperparameters given by 
\STATEx \hspace{0.5cm}  $U_{OPT} ^{(n^{(k)}-1)}$ and $\tilde \xi_{OPT}^{(n^{(k)}-1)}$  
\ENDIF
\STATE $U^{(k)}\leftarrow$ first $n^{(k)}$ singular vectors of  $ \hat A_l^{(k)} \left. \hat A_l^{(k)}\right.^\top$
\STATE $\tilde \xi^{(k)} \leftarrow \underset{\tilde \xi}{\mathrm{argmin}}\,\ell(\yv^+,\tilde \xi,U^{(k)})$ 
\REPEAT
\STATE $\tilde \xi_{OPT}^{(n^{(k)})} \leftarrow \tilde \xi^{(k)}$
\STATE $U_{OPT}^{(n^{(k)})} \leftarrow U^{(k)}$
\STATE $k \leftarrow k+1$
\STATE $n^{(k)}\leftarrow n^{(k-1)}$ 
\STATE $\hat A_l^{(k)}\leftarrow [\, \hat L_1^{(k)}\;\;  \hat L_2^{(k)}\;\; \ldots \,] $ where  $ \hat L_1^{(k)}, \hat L_2^{(k)},\ldots$ are 
\STATEx \hspace{0.5cm} the coefficients of $ \hat L^{(k)}(z)$ estimated from (\ref{Tickho_SL}) 
\STATEx \hspace{0.5cm}  with $K_S$ and $K_L$ having hyperparameters given by 
\STATEx \hspace{0.5cm}   $U^{(k-1)}$ and $\tilde \xi^{(k-1)}$  
\STATE $U^{(k)}\leftarrow$ first $n^{(k)}$ singular vectors of  $ \hat A_l^{(k)} \left.\hat  A_l^{(k)}\right.^\top$
\STATE $\tilde \xi^{(k)} \leftarrow \underset{\tilde \xi}{\mathrm{argmin}}\,\ell(\yv^+,\tilde \xi,U^{(k)})$ 
\UNTIL{ $\ell(\yv^+,\tilde \xi^{(k)},U^{(k)})<\ell(\yv^+,\tilde \xi^{(k-1)},U^{(k-1)})$ }
\UNTIL{ }
\STATEx \hspace{0.2cm} $\ell(\yv^+,\tilde \xi_{OPT}^{(n^{(k)})},U_{OPT}^{(n^{(k)})})<\ell(\yv^+, \tilde \xi_{OPT}^{(n^{(k)}-1)},U_{OPT}^{(n^{(k)}-1)})$
\STATE $n  \leftarrow n^{(k)}-1$
\STATE $U\leftarrow U_{OPT}^{(n^{(k)}-1)}$
\STATE $\tilde \xi\leftarrow  \tilde \xi^{(n^{(k)}-1)}_{OPT}$
\end{algorithmic}
\end{algorithm}
Finally, to efficiently compute a local minimum of (\ref{loglik}) in the algorithm we used the
scaled gradient projection algorithm developed in \cite{BONETTINI_2014}.

It is worth noting that our algorithm is similar to the non-separable reweighting scheme proposed in \cite{WIPF_2010} for solving a sparse Bayesian
learning problem. That algorithm iteratively alternates the computation of the optimal estimate (in our case $\hat L(z)$) and   
the closed form update of the hyperparameters (in our case it is given by (\ref{Lambda_struct})).

\section{Numerical experiments}\label{section_simulation}
We consider two Monte Carlo studies of $100$ runs where at any run a manifest process $y$ of dimension $m=6$ is considered.
For each run in the Monte Carlo experiments an identification data set of size 500 and a test set of size 1000 are generated. 
We compare the true model, denoted by TRUE, with the following estimators
\begin{itemize}
  \item SL-I: this is the sparse plus low-rank estimator (\ref{Bayes_SL}) with $K_S$ and $K_L$ as in (\ref{ME_kernel1}) 
   \item SL-II: this is the sparse plus low-rank estimator (\ref{Bayes_SL}) with $K_S$ and $K_L$ as in (\ref{ME_kernel2})
  \item L-I: this is the low-rank estimator (\ref{Bayes_SL}) with $K_S$ set equal to zero and $K_L$ as in (\ref{ME_kernel1})
  \item L-II: this is the low-rank estimator (\ref{Bayes_SL}) with $K_S$ set equal to zero and $K_L$ as in (\ref{ME_kernel2})
  \item S: this is the sparse estimator (\ref{Bayes_SL}) with $K_L$ set equal to zero
  \item SS: this is the estimator based on model (\ref{OEmodel}) where $G(z)$ is modeled as a zero-mean Gaussian process with kernel $K=I_{p^2}\otimes P$  and $P$ defined in (\ref{def_K_tilde}).
\end{itemize}
For implementation purposes, the impulse responses are truncated to a certain length $T$. The latter represents the ``practical'' length of those impulse responses \citep{PILLONETTO_2011_PREDICTION_ERROR}.
Note that, the truncation does not introduce bias-variance tradeoff because the filtered kernel $P$ forces the estimated impulse responses to decay exponentially. In these experiments we used $T=50$.

The following performance indexes are considered:
\begin{itemize}
  \item Relative complexity of the estimated model. It is quantified with 
  \al{C=\lim_{T\rightarrow \infty}\frac{\# \Ec_T }{m^2 T}}
 where $\# \Ec_T$ is the minimum number of parameters needed to characterize the  model 
 with impulse responses truncated at $T$. For instance, for the S+L model (\ref{S+Lmodel}) with $S(z)$ having $s$ nonnull transfer functions and $L(z)=FH(z)$ having rank $n$, we have $\#\Ec_T=sT+mn(T+1)$. The denominator $m^2 T$ is the number of coefficients needed in unstructured predictor model $G(z) = \sum_{k=1}^T G_k z^{-k}$ in  (\ref{OEmodel}) 

If we quantify the complexity of the S+L network as the number of edges among the manifest nodes and from the factor  ones to the manifest ones, then it is not difficult to see that its complexity is equal to $m^2 C$. Similar conclusions can be found for networks having only the sparse or the low-rank part. Therefore, the smaller $C$ is, the simpler the network  is.  
  \item One-step-ahead Coefficient of Determination, denoted by ${\mathrm{COD}_1}$. Such index quantifies how much  of the test set variance is explained by the forecast, and is defined as:
  \al{  {\mathrm{COD}_1}=1-\frac{\frac{1}{1000}\sum_{t=1}^{1000} \|y^{\mathrm{test}}(t)-\hat y^{\mathrm{test}}(t|t-1)\|^2}{\frac{1}{1000}\sum_{t=1}^{1000}\|  y^{\mathrm{test}}(t)- \bar y^{\mathrm{test}}\|^2}}
where $\bar y^{\mathrm{test}}$ denotes the sample mean of the test set data $y(1)^{\mathrm{test}}\ldots y^{\mathrm{test}}(1000)$ and $\hat y^{\mathrm{test}}(t|t-1)$ is the one-step ahead prediction computed using the estimated model. Notice that, the larger ${\mathrm{COD}}_1$ is, the better predictive performance of the estimator is;
\item Average impulse response fit
\al{  {\mathrm{AIRF}}=100 \left( 1-\sqrt{ \frac{ \frac{1}{T} \sum_{k=1}^T \| G_k -\hat G_k\|^2 }{ \frac{1}{T} \sum_{k=1}^T \| G_k -\bar G\|^2 }  }\right)}
where $G_k$ and $\hat G_k$ are the impulse response coefficients of the true and estimated model, respectively, and $\bar G=\frac{1}{T}\sum_{k=1}^T G_k$. In the case the true (estimated) model 
is S+L we have $G_k=S_k+L_k$ ($\hat G_k=\hat S_k+\hat L_k$).
\end{itemize}

 \begin{figure}[htbp]
\begin{center}
\includegraphics[width=\columnwidth]{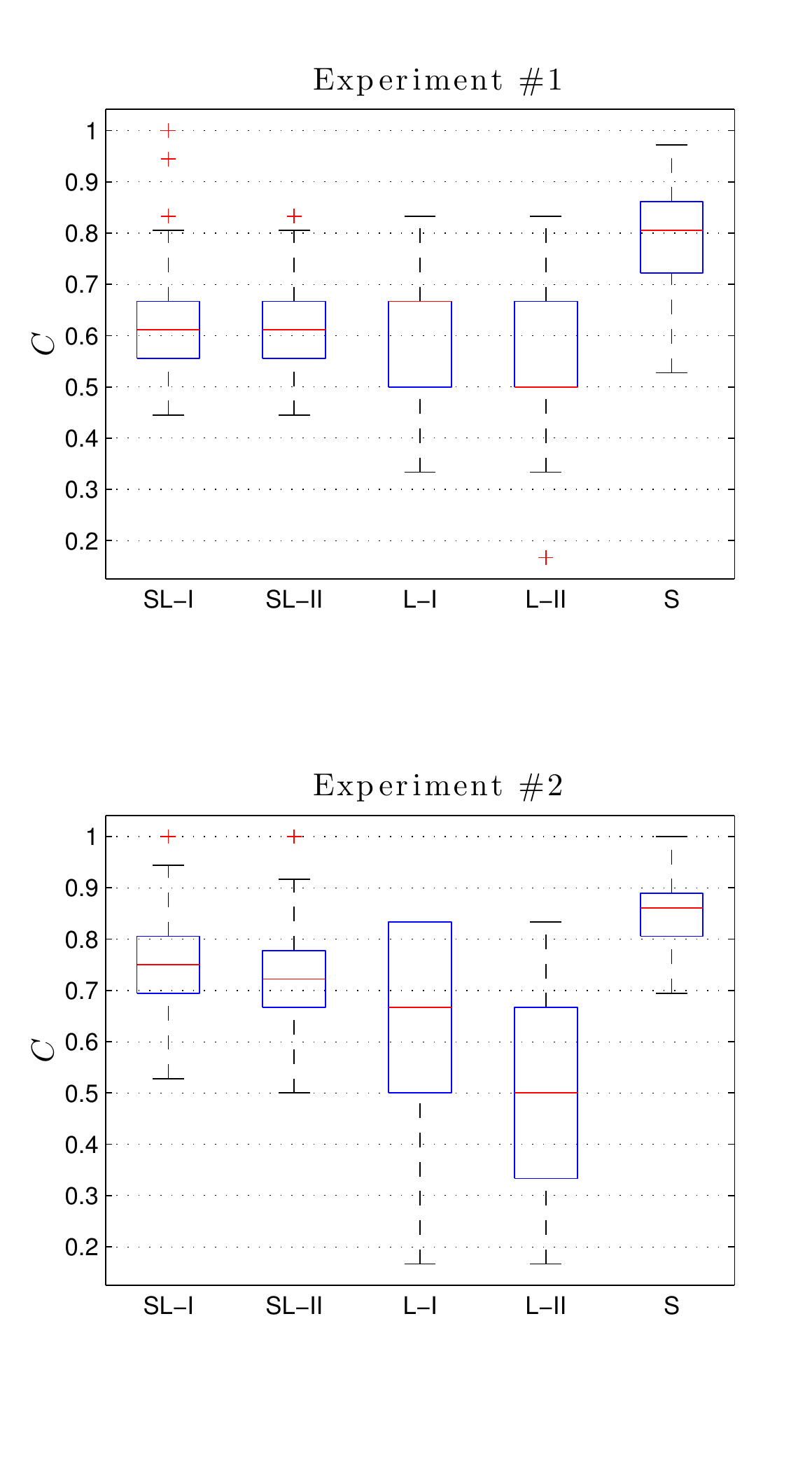}
\end{center}
\caption{Relative complexity of the models obtained by using the estimators SL-I, SL-II, L-I, L-II, and S.}\label{complex}
\end{figure}

 \begin{figure}[htbp]
 \begin{center}
\includegraphics[width=\columnwidth]{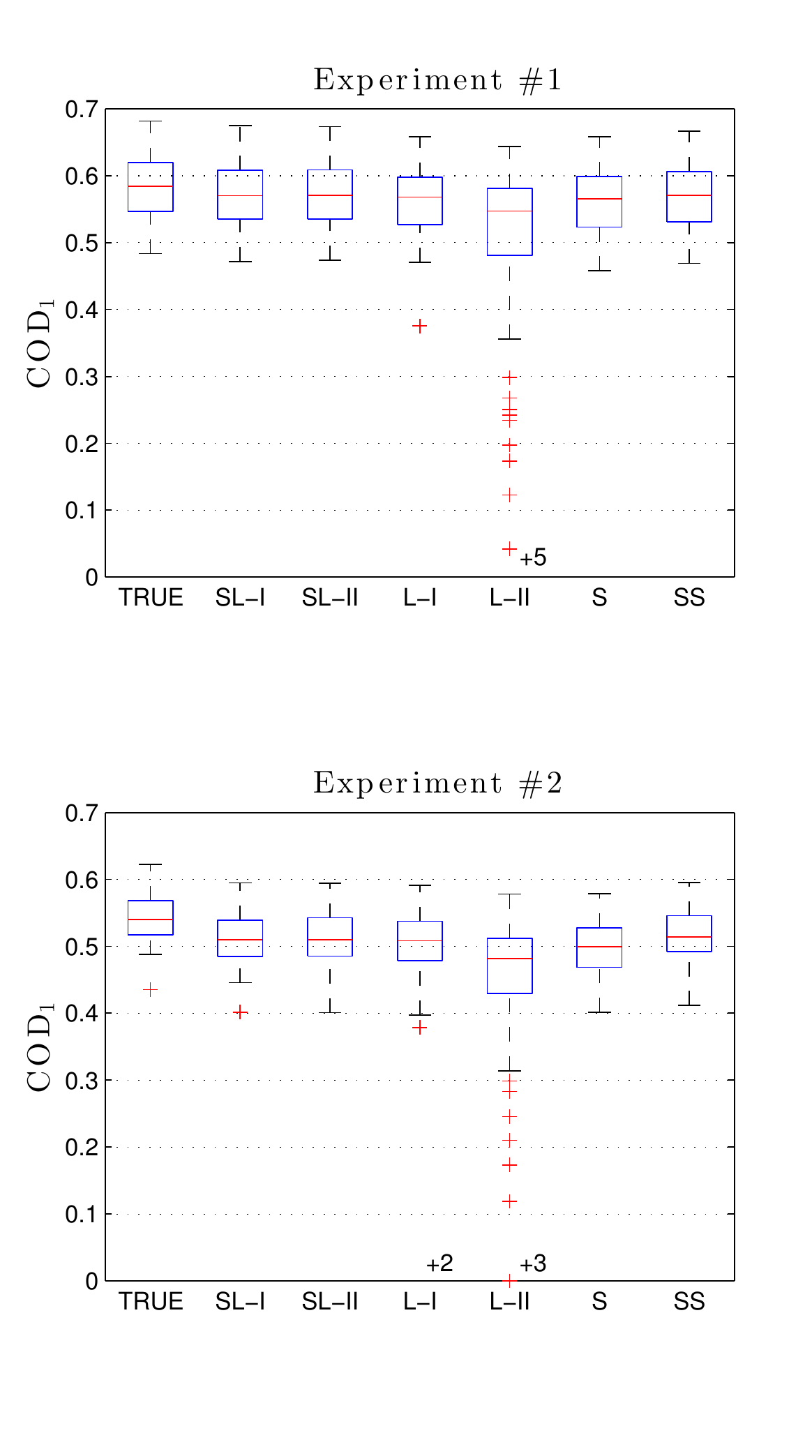}
\end{center}
\caption{One step ahead coefficient of determination ($\mathrm{COD}_1$) obtained by the 6 estimators described in Section \ref{section_simulation}. Moreover, TRUE provides an upper bound on the performance of those estimators, being the true model for the manifest process $y$.}\label{cod1}
\end{figure}

\begin{figure}[htbp]
\begin{center}
\includegraphics[width=\columnwidth]{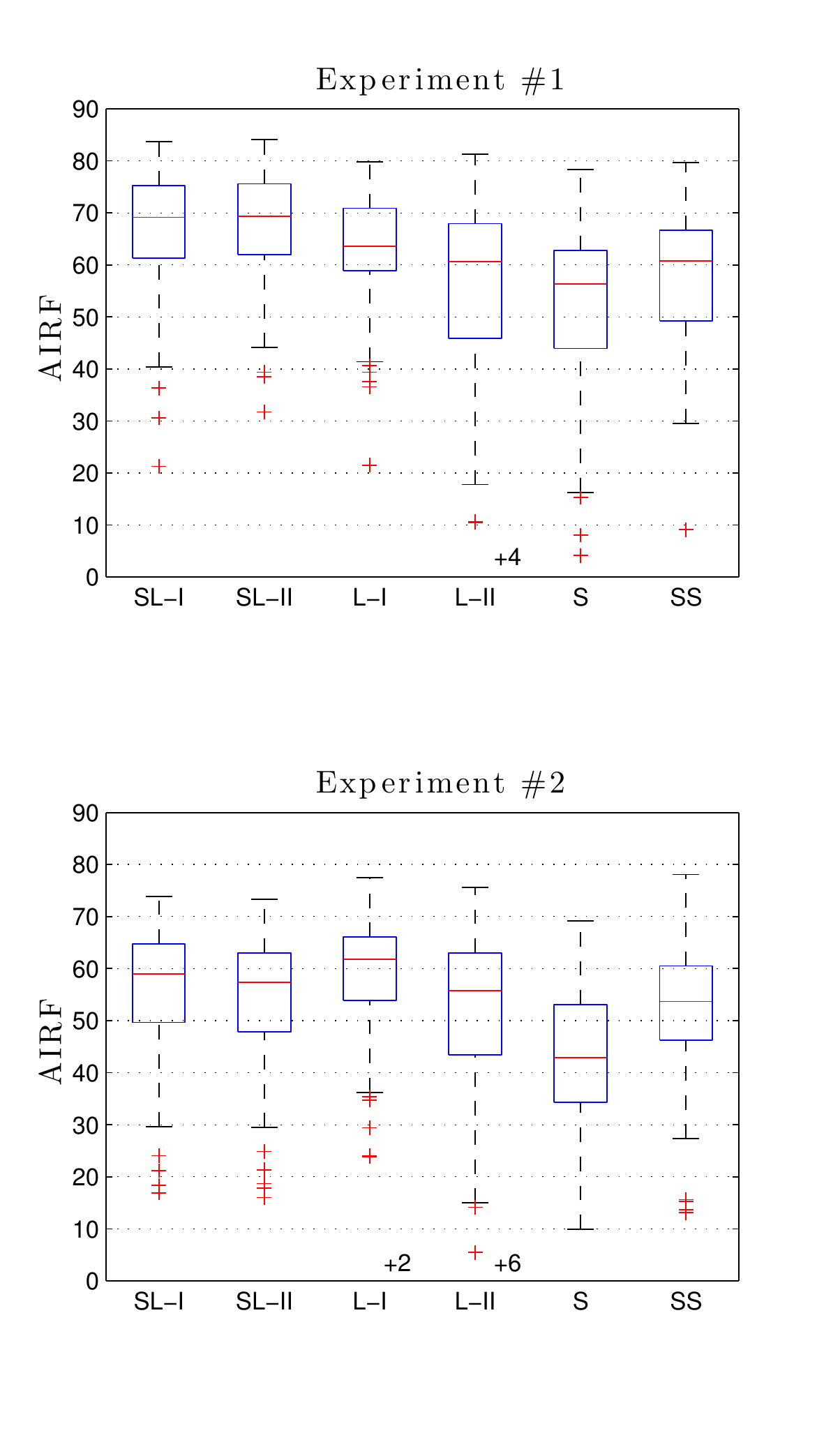}
\end{center}
\caption{Average impulse response fit ($\mathrm{AIRF}$) obtained by the 6 estimators described in Section \ref{section_simulation}.}\label{air}
\end{figure}
In the first Monte Carlo experiment, $y$ is generated through model (\ref{S+Lmodel}). In particular, $L(z)$ is randomly generated with rank $n=1$. $S(z)$ is randomly generated with 7 nonnull transfer functions. The position of those transfer functions is randomly chosen.  The top panel of Figure \ref{complex} 
shows the relative complexity of the models obtained with the estimators SL-I, SL-II, L-I, L-II and S reporting the boxplots of the values of $C$ after the 100 runs. One can see that 
L-II is the best estimator, then SL-I and S-II, then L-I and finally S. The top panel of Figure \ref{cod1}, reporting  the boxplots of  ${\mathrm{COD}}_1$, 
shows the one-step predictive performance of the estimators is similar.
The top panel of Figure \ref{air} 
shows the boxplots of $\mathrm{AIRF}$. It is clear that SL-I and SL-II outperform the others.

In the second experiment, $y$ is generated through a generic (i.e. unstructured) model.  
The bottom panels of Figure \ref{complex}, Figure \ref{cod1} and Figure \ref{air} show the performance of the employed estimators as before.  One can see that the best estimators in terms of ${\rm COD}_1$ are  SL-I, SL-II and L-I. 
SL-I and SL-II provide simpler  models  than those obtained by L-I. On the other hand in terms of  AIRF, SL-I and SL-II are slightly worse than  L-I.  
These two  Monte Carlo experiments suggest that the proposed S+L models and estimators (SL-I and SL-II) provide an effective way of estimating complex model, yielding a good tradeoff among model complexity, prediction accuracy and average impulse response fit.  

\section{Conclusions} \label{section_conclusions}
In this paper, we proposed two procedures, based on a nonparametric Gaussian regression approach, to identify stochastic processes having a sparse plus low-rank (S+L) network. The kernels inducing the S+L structure have been derived  using the maximum entropy principle. Simulations show  that
the proposed S+L estimators have  good predictive capability as well as low complexity compared with the sparse estimator and the low-rank estimators.

\section*{Appendix} 
 \subsection*{Proof of Proposition \ref{prop_theta_SL}}
  
 It is sufficient to observe that Problem (\ref{Tickho_SL}) can be rewritten as
 \al{ \label{Tikho_pb_theta}\hat \theta=\underset{\substack{\tilde \theta\in\Hc_{\tilde K}}}{\arg\min}\|\yv^+ -\tilde \Phi \tilde \theta\|^2_{\Sigma^{-1}\otimes I_N}+\| \tilde \theta\|^2_{\tilde K^{-1}}}
   where $\tilde  \theta=\left[
                                                                                                               \begin{array}{cc}
                                                                                                                 \theta_s^\top & \theta_l^\top \\
                                                                                                               \end{array}
                                                                                                             \right]^\top$, the regression matrix is  $\tilde \Phi=\left[
                                                                                                                               \begin{array}{cc}
                                                                                                                                 \Phi & \Phi \\
                                                                                                                               \end{array}
                                                                                                                             \right]$, and $\Hc_{\tilde K}$ is the reproducing Hilbert space of deterministic functions on $\Ns$ with kernel     \al{ \tilde  K=\left[
        \begin{array}{cc}
          K_S & 0 \\
          0 & K_L \\
        \end{array}
      \right].
} Then, the statement follows from Proposition 3 in  \cite{PILLONETTO_2011_PREDICTION_ERROR}. More precisely, it is not difficult to see that Proposition 3 still holds when the covariance matrix of the noise vector is not diagonal. \qed\\

 \subsection*{Proof of Proposition \ref{prop_stability}}
  Define $\breve P\in\Sc_2(\Ns)$ such that 
 \al{ [\breve{P}]_{ts}=  \left\{\begin{array}{cc} \kappa t^{\alpha} e^{-\beta t} & t=s  \\  
\left. [ P ]\right. _{ts} & \hbox{otherwise,}\\\end{array}\right.} 
hence $ P\leq  \breve P$. Let $\varepsilon>1$ and consider the decomposition 
\al{\breve P=D\bar P D}
where $D\in\Sc_2(\Ns)$ is an infinite diagonal matrix with $\sqrt{t^{\alpha+\varepsilon}}e^{-\frac{\beta}{2}t}$ as $t$-th element in the main diagonal, and 
\al{ [\bar{P}]_{ts}=  \left\{\begin{array}{cc} \kappa t^{-\varepsilon} & t=s  \\  
\star & \hbox{otherwise.}\\\end{array}\right. }Note that, $\bar P$ is strictly positive definite by construction. Moreover,
\al{\tr(\bar P)=\kappa \sum_{t =1}^\infty t^{-\varepsilon} <\infty}
because $\varepsilon>1$. 
This means the sum of all the (nonnegative) eigenvalues of $\bar P$ is bounded, therefore its maximum eigenvalue, say $\lambda$, is bounded. Therefore, by (\ref{condps}) we have 
\al{ c &\geq \Es[\phi^{\top} P^{-1}\phi]\geq \Es[\phi^{\top} \breve P^{-1}\phi]\nn\\
& \geq \Es[\phi^{\top} D^{-1}\bar P^{-1} D^{-1}\phi]\nn\\
& \geq \lambda^{-1} \Es[\phi^{\top} D^{-2}\phi]=\lambda^{-1}\sum_{t=1}^\infty t^{-(\alpha+\varepsilon)}e^{\beta t} [K]_{tt}
}  accordingly $\sum_{t=1}^\infty t^{-(\alpha+\varepsilon)}e^{\beta t} [K]_{tt}$ is bounded and thus $[K]_{tt} \leq \bar \kappa_\varepsilon t^{\alpha+\varepsilon}e^{-\beta t}$, $t\in\Ns$, for some $\bar \kappa_\varepsilon>0$.
\qed\\

  \subsection*{Proof of Theorem \ref{teo_ME}}
  
 We characterize the optimal solution by exploiting  duality theory. We consider the Lagrange function
\al{ \Lc(\pp_{\Ic},\tilde \Gamma,\lambda,\tilde \Lambda)=&\mathbf{H}(\pp_{\Ic}) \nn\\
& +\frac{1}{2}\sum_{i,j=1}^p \tilde \gamma_{(i-1)p+j} (c_{ij}- \Es[\|\check s^{[ij]}\|^2_{P^{-1}_{\Ic_s} }] )\nn\\
& +\frac{1}{2}\lambda  (r- \sum_{i,j=1}^p \Es[\|\check l^{[ij]}\|^2_{P^{-1}_{\Ic_l} }] )\nn\\
&+\frac{1}{2} \tr(\tilde \Lambda (Q-\Es[\check A_l \check A_l^\top]))}
where $\tilde \Gamma=\mathrm{diag}(\tilde \gamma_{1}\ldots \tilde \gamma_{p^2})$ with $\tilde \gamma_{i} \geq 0$, $i=1\ldots p^2$, $\lambda\geq 0$ and $\tilde \Lambda\in\overline{\Mc}_+^p$ are the Lagrange multipliers.
It is not difficult to see that $\Lc$ is strictly concave over $\Pc$. Moreover, its unique maximum point is given by annihilating its first derivative. Therefore, we obtain
\al{ &\pp_{\Ic}(\check \theta_s,\check \theta_l)=\frac{1}{c} \mathrm{exp}\left( -\frac{1}{2}\sum_{i,j=1}^p \tilde \gamma_{(i-1)p+j} \| \check s^{[ij]}\|^2_{P^{-1}_{\Ic_s} } \right.\nn\\ 
& \left. - \frac{1}{2}\lambda \sum_{i,j=1}^{p} \| \check l^{[ij]} \|^2_{P_{\Ic_l}^{-1}} - \frac{1}{2}\tr(\tilde \Lambda \check A_l \check A_l^\top)\right)} 
where $c$ is the normalization constant.
Let $e_{i}$, $i=1\ldots p^2$, denote the $j$-th vector of the canonical basis of $\Rs^{p^2}$. Then, we have
\al{ \sum_{i,j=1}^p  & \tilde \gamma_{(i-1)p+j}\| \check s^{[ij]}\|^2_{P^{-1}_{\Ic_s}}\nn\\
&=\sum_{i,j=1}^p \tilde \gamma_{(i-1)p+j} \check \theta_s^\top (e_{(i-1)p+j} e_{(i-1)p+j}^\top\otimes P^{-1}_{\Ic_s}) \check\theta_{s}\nn\\
&=\check \theta_s^\top (\tilde \Gamma\otimes P^{-1}_{\Ic_s})\check \theta_s.}
In similar way,
\al{\sum_{i,j=1}^{p} \| \check l^{[ij]} \|^2_{P_{\Ic_l}^{-1}} =\check \theta_l^\top (I_{p^2}\otimes P_{\Ic_l}^{-1}) \check \theta_l.}

Moreover, it is not difficult to see that
\al{ \label{link_thetal_Al}\tr(\tilde \Lambda \check A_l \check A_l^\top)=\check \theta_l^\top(\tilde \Lambda \otimes I_{p |\Ic_l|}) \check \theta_l. } Therefore, the optimal solution (if it does exist) is such that $\pp_{\Ic}=\pp_{\Ic_s}\pp_{\Ic_l}$, where \al{ \pp_{\Ic_s}(\check \theta_s)&=\frac{1}{c_s}\mathrm{exp}\left(-\frac{1}{2} \check \theta_s^\top (\tilde \Gamma \otimes P^{-1}_{\Ic_s})\check \theta_s \right)\nn\\
\pp_{\Ic_l}(\check \theta_l)&=\frac{1}{c_l}\mathrm{exp}\left(-\frac{1}{2}\check \theta_l^\top(\lambda I_{p^2}\otimes P_{\Ic_l}^{-1}+ \tilde \Lambda \otimes I_{p |\Ic_l |}) \check \theta_l\right)\nn }
with $c_s$ and $c_l$ normalization constants. 
Note that, $\pp_{\Ic_s}$ and $\pp_{\Ic_l}$ denote the marginal probability density of $\check \theta_s$ and $\check \theta_l$, respectively,
Therefore, $\check \theta_s$ and $\check \theta_l$ are independent, Gaussian  with zero mean and covariance matrix
\al{ \check K_S&= \tilde \Gamma ^{-1}\otimes P_{\Ic_s}\nn\\
 \check K_L&= (\lambda I_{p^2}\otimes P_{\Ic_l}^{-1}+\tilde \Lambda \otimes I_{p |\Ic_l |})^{-1} .}
Next, we prove the existence of such a solution showing that the dual problem does admits solution. Note that,
\al{ \mathbf{H}(\pp_{\Ic})=\frac{1}{2}\log\det(\check K_S)+\frac{1}{2}\log\det(\check K_L)+ \hbox{ const. term}.\nn} Therefore, the dual problem is equivalent to minimize the function
\al{ \label{dual_function}J( & \tilde \Gamma,\lambda,\tilde \Lambda)= -|\Ic_s | \log\det(\tilde \Gamma)+\tr(\tilde \Gamma C)+\lambda r\nn\\ & +\tr(\tilde \Lambda Q) -\log\det (\lambda I_{p^2}\otimes P_{\Ic_l}^{-1}+\tilde \Lambda \otimes I_{p|\Ic_l|}) } 
where  \al{\label{def_C} C=\mathrm{diag}(c_{11},\ldots,c_{1p},\ldots,c_{p1},\ldots,c_{pp}).}
Since $C\in\Mc_+^{p^2}$, $r>0$ and $Q\in\Mc_+^p$, it is not difficult to see that $J$ is lower bounded. It takes 
infinite value if and only if $\tilde \Gamma$ and/or $\lambda$ and/or $\tilde \Lambda$ are not bounded (the formal proof follows the one of Proposition 5.1 in \cite{OTTSIGMA}, see also \cite{BETA} and \cite{BETAPRED}). Moreover, if $\tilde \Gamma$ and/or $\lambda I_{p^2}\otimes P^{-1}_{ |\Ic_l |} +\tilde \Lambda\otimes I_{p |\Ic_l |}$ tend to be singular then $J$ approaches infinity. Accordingly, we can restrict the search of the Lagrange multipliers over the closed and bounded set
\al{\{ (\tilde \Gamma,\lambda,\tilde \Lambda) \hbox{ s.t. } & \varepsilon_1 I_{p^2}\leq  \tilde \Gamma\leq M_1 I_{p^2},\nn\\
&   \varepsilon_2 \leq  \lambda \leq M_2,	\;\;  \varepsilon_3 I_{p}\leq  \tilde \Lambda \leq M_3 I_{p}  \}}
for some and $\varepsilon_1,M_1>0$ and $\varepsilon_2,\varepsilon_3,M_2,M_3\geq 0$. The latter is a compact set because we are in a finite dimensional space.     
Since $J$ is continuous over that set, by the Weiestrass Theorem, the dual problem admits solution. Accordingly, the kernel matrices for $\check \theta_s$ and $\check \theta_l$
solution to (\ref{ME_problem}) does exist and are unique. Finally, setting $\gamma_i=\tilde \gamma_i^{-1}$ we obtain matrix $K_S$ in the statement. \qed\\

 \subsection*{Proof of Proposition \ref{extended_ME} }

Consider the maximum entropy problem (\ref{ME_problem}) where $c_{ij}$ and $Q$ have been replaced by $c_{ij}^{(k)}$ and $Q^{(k)}$ arbitrarily extracted from some sequences $\{c_{ij}^{(k)}\}_{k\geq 0}$, $c_{ij}^{(k)}>0$, 
and $\{Q^{(k)}\}_{k\geq 0}$, $Q^{(k)}\in\Mc^p_+$, respectively, and such that $c_{ij}^{(k)}\rightarrow 0$ $\forall\, (i,j)\in \mathbf{C}$ and $Q^{(k)}v\rightarrow 0$
$\forall\, v\in\mathbf{Q}$ as $k\rightarrow \infty$. Let $[\check \theta_s^{(k)\top}\; \check \theta_l^{(k)\top}]^\top$ be the random vector solution to this maximum entropy problem. The corresponding kernel matrices are (\ref{optimal_Ks}) and (\ref{optimal_Kl}) where $\gamma_i$, $\lambda$ and $\Lambda$
have been replaced with $\gamma^{(k)}_i>0$, $\lambda^{(k)}>0 $ and $\Lambda^{(k)}\in \Mc^p_+$ with $k\geq 0$.
Then, substituting this random vector into (\ref{constr_s}) we obtain
\al{c_{ij}^{(k)} & \geq \Es[\| \check s^{[ij]}\|^2_{P_{\Ic_l}^{-1}}]=\tr(\Es[ \check s^{[ij]} \check s^{[ij]^\top}] P_{\Ic_s}^{-1})\nn \\ &=|\Ic_s| \gamma_{(i-1)p+j}^{(k)}> 0.}
Accordingly, if $(i,j)\in\mathbf{C}$ then $\gamma^{(k)}_{(i-1)p+j}\rightarrow 0$ as $k\rightarrow \infty$. Let $\check A_l^{(k)}$ be the matrix built from the maximum entropy random vector.  Then, substituting it in constraint (\ref{constr_l2}) pre- and post-multiplied by an arbitrary $v\in\Rs^{p}$, we have   
 \al{v^\top & Q^{(k)} v\geq v^\top \Es[\check A_l^{(k)} \check A_l^{(k)\top} ]v=\Es[\tr(v v^\top \check A_l^{(k)} \check A_l^{(k)\top} )]\nn\\ 
&= \Es[\check \theta_l^{(k)\top}(vv^\top\otimes I_{p|\Ic_l|}) \check \theta_l^{(k)}]\nn\\ &=
\tr((vv^\top\otimes I_{p|\Ic_l|})\Es[\check \theta_l^{(k)}\check \theta_l^{(k)\top}] ) \nn\\
& =\tr((v v^\top \otimes I_{p|\Ic_l |})(\lambda^{(k)} I_{p^2}\otimes  P_{\Ic_l}^{-1}+\Lambda^{(k)^{-1}} \otimes I_{p |\Ic_l |} )^{-1})\nn\\
&\geq v^\top (\lambda^{(k)}\mu I_p+ \Lambda^{(k)^{-1}})^{-1} v p|\Ic_l| \nn\\
&= v^\top ( \Lambda^{(k)}-    \Lambda^{(k)}( \Lambda^{(k)}+\frac{1}{\lambda^{(k)}\mu }I_p)^{-1}  \Lambda^{(k)} )   v p|\Ic_l| > 0}
where we exploited (\ref{link_thetal_Al}) and $\mu>0$ is the maximum eigenvalue of $P_{\Ic_l}$. Accordingly, if $v\in\mathbf{Q}$ then $\Lambda^{(k)}v\rightarrow 0$ as $k\rightarrow \infty$. \qed\\

 \subsection*{Proof of Proposition \ref{prop_stability2}}

 We use the same decomposition for $\breve P$ exploited in the proof of Proposition \ref{prop_stability}. Thus, we have
\al{ Q& \geq \Es[A_l (P^{-1}\otimes I_p)A_l^\top] \geq \Es[A_l (\breve P^{-1}\otimes I_p)A_l^\top]\nn\\ 
&\geq \Es[A_l (D^{-1}\otimes I_p)(\bar P^{-1}\otimes I_p)(D^{-1}\otimes I_p)A_l^\top]\nn\\
&\geq \lambda^{-1}\Es[A_l (D^{-2}\otimes I_p)A_l^\top]\nn\\
&\geq \lambda^{-1}\sum_{t=1}^\infty t^{-(\alpha+\varepsilon)}e^{\beta t} \Es[L_t L_t^\top]}
where $\lambda$ is the maximum eigenvalue of $\bar P$ which is bounded. Hence, condition (\ref{constraint_p_l3}) implies that
\al{\sum_{t=1}^\infty t^{-(\alpha+\varepsilon)}e^{\beta t} \Es[|[L_t]_{ij}|^2]}
is bounded for $i,j=1\ldots p$. Accordingly, $\Es[| [L_t]_{ij}|^2]\leq \bar\kappa_\varepsilon t^{\alpha+\varepsilon}e^{-\beta t}$, $t\in\Ns$, for some $\bar \kappa_\varepsilon>0$. 
Last, let $v$ be in the null space of $Q$, i.e. $v^\top Q v=0$. It follows from \eqref{constraint_p_l3} that 
$$
\Es[v^\top A_l  (P^{-1}\otimes I_p) A_l^\top v] = v^\top Q v = 0
$$
Since $(P^{-1}\otimes I_p) $ is positive definite this implies that $\Es[v^\top A_l A_l^\top v] = 0$, which completes the proof.\qed\\

 \subsection*{Proof of Theorem \ref{teo_ME2}}
The statement is proved by using the duality theory as in the proof of Theorem \ref{teo_ME}. In particular the Lagrange function is 
\al{\Lc(\pp,\tilde \Gamma,\tilde \Lambda)=&\mathbf{H}(\pp_{\Ic}) \nn\\
& +\frac{1}{2}\sum_{i,j=1}^p \tilde \gamma_{(i-1)p+j} (c_{ij}- \Es[\|\check s^{[ij]}\|^2_{ P^{-1}_{\Ic_s} }] )\nn\\
& +\frac{1}{2} \tr(\tilde \Lambda(Q-\Es[\check A_l( P^{-1}_{\Ic_l} \otimes I_p)\check A_l^\top]))	 }
where $\tilde \Gamma =\mathrm{diag}(\tilde \gamma_1\ldots \tilde \gamma_{p^2})$ with $\tilde \gamma_i\geq 0$, $i=1\ldots p^2$, and $\tilde \Lambda\in\overline{\Mc}_p^+$. 
Then, it is not difficult to see that $\pp_{\Ic}$ maximizing $\Lc$ is such that $\pp_{\Ic}=\pp_{\Ic_s}\pp_{\Ic_l}$
where \al{\pp_s&=\frac{1}{c_s}\mathrm{exp}\left(-\frac{1}{2}\check \theta_s^\top (\tilde \Gamma\otimes  P^{-1}_{\Ic_s})\check\theta_s \right)\nn\\
\pp_l&=\frac{1}{c_l}\mathrm{exp}\left(-\frac{1}{2}\check \theta_l^\top(\tilde \Lambda\otimes I_p \otimes P^{-1}_{\Ic_l}) \check \theta_l\right)}
with $c_s$ and $c_l$ normalization constants. Therefore, the optimal solution (if it does exist) is such that $\check \theta_s$ and $\check \theta_l$ are independent, Gaussian, with zero mean and covariance matrix, respectively,  
\al{K_S&= \tilde \Gamma^{-1}\otimes P_{\Ic_s}\nn\\
K_L& = \tilde \Lambda^{-1}\otimes I_p \otimes P_{\Ic_l}.}
The existence of such a solution is proved by showing that the dual problem admits solution. The latter consists in minimizing the function 
\al{J(\tilde \Gamma,\tilde \Lambda)= - |\Ic_s |&  \log\det(\tilde \Gamma)+\tr(\tilde \Gamma C)\nn\\ &  -p|\Ic_l| \log\det (\tilde \Lambda) +\tr(\tilde \Lambda Q)} 
where $C$ has been defined in (\ref{def_C}).
In this case the search of the minimum can be restricted to the compact set 
\al{ \{(\tilde \Gamma,\tilde \Lambda) \hbox{ s.t. } 0<\tilde \Gamma \leq M_1 I_{p^2},\;\; 0<\tilde \Lambda \leq M_2 I_p\}}
for some $M_1,M_2>0$. Moreover, $J$ is continuous over this set. Thus, by the Weiestrass Theorem the dual problem admits solution. 
Finally, setting $\gamma_i=\tilde \gamma_i^{-1}$ and $\Lambda=\tilde \Lambda^{-1}$  we obtain the kernel matrices in the statement. \qed\\

\subsection*{Proof of Proposition \ref{extended_ME2}}
The proof follows the same lines of that of Proposition \ref{extended_ME}.


\end{document}